\documentclass[10pt]{article}

\usepackage[unicode]{hyperref}
\usepackage{amsmath}
\usepackage{graphicx}								
\usepackage{subfig}
\usepackage{amsfonts} 							
\usepackage{dsfont}
\usepackage{algorithm}              
\usepackage{algorithmic}
\usepackage{booktabs} 							
\usepackage{multirow}								
\usepackage{setspace}
\usepackage{mathrsfs}
\usepackage{amssymb}
\usepackage{tikz}
\usepackage{color}
\usepackage{xcolor}
\usepackage{mathtools}
\usepackage{xfrac}
\usepackage{upgreek}

\usepackage{framed}
\usepackage{afterpage}
\usepackage{capt-of}

\usepackage{geometry}
\usepackage{color}
\definecolor{deepblue}{rgb}{0,0,0.5}
\definecolor{deepred}{rgb}{0.6,0,0}
\definecolor{deepgreen}{rgb}{0,0.5,0}
\definecolor{gray(x11gray)}{rgb}{0.75, 0.75, 0.75}

\usepackage{tikz}	
\usetikzlibrary{arrows,chains,matrix,positioning,scopes}

\makeatletter
\tikzset{join/.code=\tikzset{after node path={%
\ifx\tikzchainprevious\pgfutil@empty\else(\tikzchainprevious)%
edge[every join]#1(\tikzchaincurrent)\fi}}}

\makeatother
\tikzset{>=stealth',every on chain/.append style={join}, every join/.style={->}}
\tikzstyle{labeled}=[execute at begin node=$\scriptstyle,   execute at end node=$]

\tikzset{
>=stealth',
help lines/.style={dashed, thick},
axis/.style={<->},
important line/.style={thick},
connection/.style={thick, dotted},
}

\newcommand {\Sum}   {\sum\limits}

\newcommand {\R}		 {\mathbf{r}}

\newcommand {\Rd}    {{\mathds{R}}^d}

\def \EnergyNorm#1  {{\mid\!\mid\!\mid #1 \mid\!\mid\!\mid}^2 }   

\def \dvrg       {\mathrm{div}}	

\def \traspose#1 {{#1}^{rm T}}
\def \laplace    {\Delta}

\newcommand {\vectorn}  {\boldsymbol {n}}

\newcommand {\flux}     {\boldsymbol {y}}

\newcommand {\hhat} {\hat{h}} 
\newcommand {\Khat} {\widehat{K}} 
\newcommand {\Bhat} {\widehat{B}} 
\newcommand {\Rhat} {\widehat{R}} 
\newcommand {\Qhat} {\widehat{Q}} 

\def \dxt      {\mathrm{\:d}x\mathrm{d}t}

\def\L#1{L^{#1}}
\def\H#1{H^{#1}}

\def\HD#1#2{H^{#1}_{#2}}

\newcommand {\CFriedrichs} {C_{{\rm F}}}

\newcommand*\rfrac[2]{{}^{#1}\!/_{#2}}

\newcommand{\mnote}[1]{\!^\textrm{\scriptsize\color{gray(x11gray)}#1}}

\newtheorem{theorem}{Theorem}{\bf}{\it}
\newtheorem{corollary}{Corollary}{\bf}{\it}
{\bf}{\it}
\newtheorem{lemma}{Lemma}{\bf}{\it}
{\bf}{}
\newtheorem{remark}{Remark}{\bf}{\it}

\def\ProofBegin{\noindent{\bf Proof:} \:}
\def\ProofEnd{{\hfill $\square$}}

\definecolor{formalshade}{rgb}{0.95,0.95,1}

\definecolor{green}{rgb}{0.09, 0.45, 0.27} 
\definecolor{magenta}{rgb}{0.5, 0.0, 0.5}
\definecolor{brown}{rgb}{0.6,0.2,0}

\title{A posteriori error estimates for space-time IgA approximations to 
parabolic initial boundary value  problems}
\author{Ulrich Langer \thanks{RICAM Linz, Johann Radon Institute, Linz, Austria} \and
Svetlana Matculevich \thanks{RICAM Linz, Johann Radon Institute, Linz, Austria} \and 
Sergey Repin \thanks{St. Petersburg Department of V.A. Steklov Institute of Mathematics RAS; 
University of Jyvaskyla, Finland}}
	
\begin{document}
	
\maketitle 

	
\begin{abstract}	
The paper is concerned with stabilised space-time approximations for initial boundary value problems
of the parabolic type. These approximations have been presented and studied by 
Langer, Neum\"uller, and Moore (2016), who have shown that they satisfy standard a priori error 
estimates. 
The goal of this paper is to deduce a posteriori error estimates and investigate their applicability to 
IgA space-time approximations.
%
%
The derivation is based on purely 
functional arguments and, therefore, the estimates do not contain mesh dependent constants and are 
valid for any approximation from the admissible (energy) class. In particular, they imply estimates for 
discrete norms associated with IgA approximations. We establish different forms of a posteriori error 
majorants and prove equivalence of them to the energy error norm.
%
%
This property justifies efficiency and reliability of a posteriori error estimates. Another important property 
of the estimates is their flexibility with respect to several free parameters. Using these parameters, we can obtain estimates for different 
error norms and minimise the respective majorant in order to find the best possible bound of the error.
 


\end{abstract}

\section{Introduction}

Time-dependent systems governed by parabolic equations arise in various scientific and engineering 
applications, 
for instance, processes with slow evolution such as heat conduction and diffusion, changing in time processes in 
life and social sciences, etc. Analytic and numerical treatment of this class of problems involve 
several complications due to possible non-linearities of the studied processes, instabilities of numerical 
schemes 
causing blow-ups in simulations, and increasing an amount of data due to extra dimensionality.
Therefore, remaining open questions (from theoretical as well as applied viewpoints) have triggered 
active investigations of such models in mathematical and numerical modelling and large-scale scientific 
computing.

Time-stepping methods, which have become quite popular in industrial software packages, allow 
combining
various discretisation techniques in space (e.g., finite element method, finite difference method, 
finite volume method) 
with marching in time. Usually, it is common to distinguish between two different time-stepping methods, 
i.e., horizontal and vertical methods of
lines. In the first one, also known as Rothe's method, one starts from discretisation w.r.t.  
time variable \cite{LMR:Lang:2001}, whereas the vertical method performs first the discretisation in 
space and then in time \cite{LMR:Thomee:2006}. For both of these approaches, development of 
efficient and fully adaptive schemes becomes complicated because of the separation of 
time and space discretisations. It also affects negatively the parallelisation of the solver due to the
curse of sequentiality.

Due to the fast development of parallel computers with hundreds of thousands of cores, 
treating time as yet another dimension in space in the evolutionary equations became quite natural. 
Moreover, space-time approach does not have the above-mentioned drawbacks of time-marching schemes. 
On the contrary, it becomes quite advantageous when efficient parallel methods and their implementation 
on massively parallel computers are concerned. The simplest ideas for space-time solvers are based on 
time-parallel integration techniques for ordinary differential equations (the comprehensive 
overview on the history of this approach can be found in \cite{LMR:Gander:2015}). Time-parallel
multigrid methods for parabolic problems have also a long history starting from the first introduction in 
\cite{LMR:Hackbusch:1984}. Later, parallel multigrid waveform relaxation methods for parabolic initial
boundary value problem (I-BVPs) was presented in \cite{LMR:LubichOstermann:1987}. 
The study on convergence behaviour of these time-parallel multigrid methods by means of Fourier mode 
analysis is presented in \cite{LMR:VandewalleHorton:1995}, and the paper 
\cite{LMR:Deshpandeetall:1995} provides a rigorous analysis of time domain parallelism. 
In \cite{LMR:GanderNeumuller:2014, LMR:GanderNeumuller:2016}, authors use Fourier analysis to 
construct perfectly scaling parallel space-time multigrid 
methods for solving initial value problems for ordinary differential equations, parabolic I-BVPs, and 
linear system arising from a space-time discontinuous Galerkin discretisations (see also related publications
\cite{LMR:NeumullerSteinbach:2013, LMR:NeumullerSteinbach:2011, LMR:KarabelasNeumuller:2015}). 
For a more detailed overview of the existing works on space-time techniques, we refer the reader to 
\cite{LMR:LangerMooreNeumueller:2016a}.

Various approximation methods have been recently developed for space-time formulation of I-BVPs. 
In particular, $h\!-\!p$ versions of the finite element method in space with $p$ and $h\!-\!p$ 
approximations in time for parabolic I-BVPs have been originally presented in 
\cite{LMR:BabuskaJanik:1989} and \cite{LMR:BabuskaJanik:1990}, respectively. Wavelet methodology 
was extended to space-time adaptive schemes in \cite{LMR:SchwabStevenson:2009}.
Uniform stability of abstract Petrov-Galerkin discretisations of boundedly invertible operators and their applicability to space-time discretizations of linear parabolic problems is discussed in 
\cite{LMR:Mollet:2014}. Error bounds for reduced basis approximation to linear parabolic problems were 
proved in \cite{LMR:UrbanPatera:2014}, and in \cite{LMR:Steinbach:2015} conforming space-time finite element approximations for the same class of problems was investigated. 

Increasing popularity of space-time methods 
has generated new methods of solving complicated engineering problems such as
%
fluid-structure interaction, aerodynamics problems, and cardiac electro-mechanics 
(see \cite{LMR:TakizawaTezduyar:2011, LMR:Takizawaetall:2012, LMR:Karabelas:2015}, 
and the references therein). These results only confirm the great potential of space-time methods for solving 
time-dependent problems, models with growing and shrinking in time domains, and objects with 
moving boundaries or interfaces. 

In this paper, we use the method presented in \cite{LMR:LangerMooreNeumueller:2016a}, 
based on special time-upwind test functions motivated by the space-time streamline diffusion method introduced in \cite{LMR:Hansbo:1994, LMR:Johnson:1987, LMR:JohnsonSaranen:1986}.
IgA framework provides approximations of high accuracy and flexibility due to a high smoothness of 
the respective basis functions (B-splines, NURBS, or localised splines; see, e.g., 
\cite{LMR:TakizawaTezduyar:2014, LMR:TakizawaTezduyar:2011}). Therefore, the method presented in \cite{LMR:LangerMooreNeumueller:2016a}, combining full space-time approach with IgA technologies, 
is a pioneering work into the direction of efficient fully-adaptive and heavily parallelised schemes aiming 
to tackle problems oriented to industrial applications. 

Investigation of effective adaptive refinement methods is highly important for the construction of fast 
and efficient solvers for partial differential equations (PDEs). In space-time methods (as in many others), 
the aspect of scheme localisation is strongly linked with reliable and quantitatively efficient a posteriori 
error estimation (the general overview on error estimators can be found in, e.g., monographs \cite{AinsworthOden2000, Malietall2014}). 
In other words, adaptive algorithms rely on a posteriori error estimation tools, which suppose to identify 
those areas of the considered computational domain, where the approximation error is substantially 
higher than on  the rest of it. A smart combination of solvers and error indicators
makes the refinement step fully automated to the characteristics of the problem 
(external forces, geometries, etc.) providing at the same time discretisations with the desired accuracy 
in terms of the output quantity of the interest. Moreover, such automation becomes quite essential in the 
generation of the mesh suitable to a complicated geometry, by using an efficient refinement  procedure 
and adapting the initial design representations.

Due to a tensor-product setting of IgA splines, mesh refinement has global effects, including a large 
percentage of superfluous control points in data analysis, unwanted ripples on the surface, etc. Arising 
from these challenges with the design process as well as complications in handling big amount of 
corresponding data naturally have triggered the development of local refinement strategies for IgA. 
There are at least three different approaches to achieve local refinements. The first one, 
so-called {\em truncated B-splines} (T-splines), was introduced in \cite{LMR:Sederbergetall2003, 
LMR:Sederbergetall2004} and analysed in \cite{LMR:Bazilevsetall2010, LMR:Veigaetall2011, 
LMR:Scottetall2011, LMR:Scottetall2012}. It is based on T-mesh that allows 
eliminating the redundant control points from NURBS model. The study of this 
approach has confirmed to generate efficient local refinement algorithm for analysis-suitable T-splines, 
which avoids excessive propagation of control points. 
%
An alternative approach, that allows local control of the refinement, is based on {\em truncated 
hierarchical B-splines} (THB-splines). The procedure to construct a basis of the hierarchical spline space 
was suggested in \cite{LMR:Kraft1997} and extended in \cite{LMR:Vuongetall2011, 
LMR:GiannelliJuttlerSpeleers2012}. Unlike T-spline localisation algorithm that does not eliminate the 
unwanted propagation of the refinement, no such propagation is observed for THB-splines (the 
corresponding examples can be found in \cite{LMR:NguyenThanhetall:2011, 
LMR:KleissJuttlerZulehner:2012}). 
The third group of locally defined splines is called {\em locally refined splines} (LR-splines) and have been 
developed in \cite{LMR:DokkenLychePettersen2013} and \cite{LMR:Bressan2013}. 

Local refinement techniques in IgA have been combined with various a posteriori error estimation 
approaches. For instance, a posteriori error estimates using the hierarchical bases (i.e., saturation 
assumption on the enlarged underlying space and the constants in the strengthened Cauchy inequality) 
was investigated in \cite{LMR:DorfelJuttlerSimeon2010, LMR:Vuongetall2011}. However, according to the 
original paper on this method \cite{LMR:BankSmith:1993}, the validity of this assumption depends 
strictly on the considered example. Moreover, an accurate estimation of constants in the strengthened
Cauchy inequality requires the solution of generalised minimum eigenvalue problem, which might become
quite technical.
%
In \cite{LMR:Johannessen:2009, LMR:Wangetall2011, LMR:BuffaGiannelli2015}, and 
\cite{LMR:KumarKvamsdalJohannessen2015}, residual-based a posteriori error estimators 
and their modifications were exploited 
in order to construct efficient automated mesh refinement algorithms. 
These estimates require computation of constants related to Clement-type interpolation operators, which are 
mesh dependent. Finally, goal-oriented 
error estimators based on auxiliary global refinement steps have been considered in 
\cite{LMR:ZeeVerhoosel2011, LMR:DedeSantos2012, LMR:Kuru:2013, LMR:Kuruetall2014}.
Below, we use a different (functional) method that provides fully guaranteed error estimates in the various weighted norms equivalent to the global energy norm. The estimates include only global constants (independent of the mesh characteristic $h$) and 
are valid for any approximation from the admissible functional space.

Functional error estimates (so-called majorants and minorants) were originally introduced in 
\cite{LMR:Repin:1997, LMR:Repin:1999, LMR:Repin:2000, LMR:Repin:2002} 
and later applied to different mathematical models (see the monographs \cite{LMR:RepinDeGruyterMonograph:2008, Malietall2014}). They provide guaranteed, sharp, and fully computable upper and lower bounds of 
errors. This approach to error control 
was applied to IgA schemes in \cite{LMR:KleissTomar2015}, where it was confirmed that the majorants 
also provide not only reliable and efficient upper bounds of the total energy error but a quantitatively sharp 
indicator of local errors.
%
%

In this paper, we deduce and study functional type a posteriori error estimates for time-dependent 
problems \cite{LMR:Repin:2002} in the context of the space-time IgA scheme
\cite{LMR:LangerMooreNeumueller:2016a}. By exploiting the universality and efficiency of the considered 
error estimates as well as taking an advantage of smoothness of the obtained approximations, one can construct
fully adaptive fast and efficient parallelised space-time method that could tackle complicated problems 
inspired by industrial applications. 

This paper is organised as follows: Section \ref{sec:model-problem} defines the problem  
and its variational formulation. It also introduces the notation and some special functional spaces used throughout the paper. An overview of main ideas and definitions used in the IgA framework can be found 
in the subsequent section. Section \ref{sec:discretization-non-moving-domain} presents the 
stabilised space-time IgA scheme and establishes its main properties. 
In Section \ref{eq:general-error-estimate}, we introduce new a posteriori error estimates on a functional 
type using the ideas coming from stabilised formulation of parabolic I-BVPs. 
Theorems  \ref{th:theorem-majorant-general-1} and \ref{th:theorem-majorant-general-2} present two 
different forms of the estimates that rely on different regularity assumptions for the approximate solution 
and auxiliary flux. 
Consequently, Corollaries \ref{cor:majorant-1} and \ref{cor:majorant-2} present majorants that are 
tailored to the space-time IgA scheme presented in Section \ref{sec:discretization-non-moving-domain}.
Finally, Section \ref{sec:equivalence} introduces the advanced form of the majorants (derived in 
Theorems \ref{th:theorem-majorant-general-1} and \ref{th:theorem-majorant-general-2}) and shows 
the equivalence of these modified functional estimates to the error measured in the energy norm.

\section{Model Problem}
\label{sec:model-problem}

Let {$\overline{Q} := Q \cup \partial Q$, $Q := \Omega \times (0, T)$}, denote the space-time cylinder, 
where  $\Omega \subset \Rd$, $d \in \{1, 2, 3\}$, is a bounded Lipschitz domain with boundary 
$\partial \Omega$, and $(0, T)$ is a given time interval, $0 < T < +\infty$. Here, the cylindrical surface 
is defined as $\partial Q := \Sigma \cup \overline{\Sigma}_{0} \cup \overline{\Sigma}_{T}$ with 
$\Sigma = \partial \Omega \times (0, T)$,  $\Sigma_{0} =  \Omega \times \{0\}$ and 
$\Sigma_{T} = \Omega \times \{T\}$. 
In order to avoid non-principal technical difficulties and present the main ideas in the most transparent 
form, we discuss our approach to guaranteed error control of space-time approximations with the 
paradigm of the classical {\em linear parabolic initial-boundary value problem}: find 
$u: \overline{Q} \rightarrow \mathds{R}$ satisfying the system
\begin{alignat}{2}
\partial_t u - \laplace_x u & = f \;\;\qquad {\rm in} \quad Q, \label{eq:equation}\\
u & = 0 \;\;\qquad {\rm on} \quad \Sigma, \label{eq:Dirichlet-boundary-condition} \\
u & = u_0 \qquad {\rm on} \quad { \overline{\Sigma}_0}, \label{eq:initial-condition}
\end{alignat}
where $\partial_t$ denotes the time derivative, 
$\laplace_x$ is the Laplace operator in space, 
$f \in \L{2}(Q)$ is a given source function, and 
$u_0 \in \HD{1}{0}(\Sigma_0)$ is a given initial state, satisfying zero boundary condition on 
$\partial \Sigma_0 = \partial \Omega \times \{ 0\}$. 
%
Here, $\L{2}(Q)$ denotes the space of square-integrable functions over $Q$. The respective 
norm and scalar product are denoted by $\| \, v \, \|_{Q} := \| \, v \, \|_{\L{2}(Q)}$ and 
$(v, w)_Q := \int_Q v(x,t) w(x,t) \dxt$, $\forall v, w \in \L{2} (Q)$, respectively, 
with similar notation used for spaces of vector-valued fields.

By $\H{k}(Q)$, $k \geq 1$, 
we denote spaces of functions having generalised square-summable derivatives of the order $k$
with respect to (w.r.t.) space and time. Next, we introduce the following Sobolev spaces
\begin{alignat}{2}
V_{0} := \HD{1}{0}(Q) & := \Big\{ \, u \in \H{1}(Q) \; : \; u = 0 \; {\rm on} \; \Sigma \,\Big\}, \nonumber\\
%
\HD{1}{0, \overline{0}}(Q) & := \Big\{ u \in \HD{1}{0}(Q)\;: \;u = 0 \; {\rm on} \; \Sigma_T\Big\}, \nonumber\\
H^{1}_{0, \underline{0}} (Q) & := \Big\{ u \in \HD{1}{0}(Q): u = 0 \; {\rm on} \; \Sigma_0 \Big\}, \nonumber\\
V^{\Delta_x}_{0} := \HD{\Delta_x, 1}{0}(Q)  & := \Big\{\, u \in \HD{1}{0}(Q) \,:\, \Delta_x u \in \L{2}(Q) \, \Big\}, 
\label{eq:h-deltax-1}
\end{alignat}
and 
\begin{equation}
H^{s}_{0, \underline{0}} (Q) := \H{s}(Q) \cap \HD{1}{0, \underline{0}}(Q).
\end{equation}
%
%
%
In addition, we introduce Hilbert spaces for auxiliary vector-valued functions (which are used in 
the derivation of the a posteriori error estimates): 
\begin{equation}
	H^{\dvrg_x, 0}(Q) := 
	\Big \{  \, \flux \in [\L{2}(Q)]^d \;:\; 
	         \dvrg_x \flux \in \L{2} (Q)
	\, \Big \}
	\label{eq:y-set-div-0}
\end{equation} 	
and 
\begin{equation}
	H^{\dvrg_x, 1}(Q) := \Big\{\, \, \flux \in [\L{2}(Q)]^d \;:\; 
	         \dvrg_x \flux \in \L{2} (Q) , \; \partial_t \flux \in [\L{2}(Q)]^d \, \Big\}.
	\label{eq:y-set-div-1}
\end{equation} 	
These spaces are supplied with the natural norm and semi-norm
$$\| \flux \|^2_{H^{\dvrg_x, 0}} := \|  \dvrg_x \, \flux \|^2_Q
\quad \mbox{and} \quad
\| \flux \|^2_{H^{\dvrg_x, 1}} := \|  \dvrg_x \, \flux \|^2_Q + \|  \partial_t \, \flux \|^2_Q.$$
%
%
In what follows,  $\CFriedrichs$ denotes the constant in the Friedrichs inequality
$$\| w \|_Q \leq \CFriedrichs \, \| \nabla_x w \|_Q, \quad 
\forall w \in \HD{1, 0}{0}(Q) := \Big\{ \, u \in \L{2}(Q) \; : \nabla_x u \in [\L{2}(Q)]^d, 
\; u = 0 \; {\rm on} \; \Sigma \,\Big\}.$$ 

%
%
It is also well-known that if $f \in \L{2}(Q)$ and $u_0 \in \HD{1}{0}(\Sigma_0)$, the problem \eqref{eq:equation}--\eqref{eq:initial-condition} is uniquely solvable in $V^{\Delta_x}_0$, and 
the solution $u$ depends continuously on $t$ in the norm $\HD{1}{0}(\Omega)$ (see, e.g., \cite{LMR:Ladyzhenskaya:1954} and \cite[Theorem 2.1]{LMR:Ladyzhenskayaetal:1967}). 
Moreover, according to \cite[Remark 2.2]{LMR:Ladyzhenskayaetal:1967}, 
$\| \, u_{x} (\cdot, t) \, \|^2_{\Omega}$ is an absolutely continuous function of $t \in [0, T]$ for any 
$u \in V^{\Delta_x}_0$. 
If the initial condition $u_0 \in \L{2}(\Sigma_0)$, then 
the problem has a unique solution $u \in \HD{1, 0}{0}(Q)$, 
that satisfies the generalised statement of the problem
\begin{equation}
a(u, w) = l(w), \quad \forall w \in \HD{1}{0, \overline{0}}(Q),
\label{eq:variational-formulation}
\end{equation}
with the bilinear form
\begin{equation}
	a(u, w) := 
	(\nabla_x {u}, \nabla_x {w})_Q - (u, \partial_t w)_Q,
	\label{eq:bilinear-form}
\end{equation}
and the linear functional
\begin{equation}
{  
l(w) := (f, {w})_Q + (u_0, {w})_{\Sigma_0}.
}	
\label{eq:linear-functional}
\end{equation}
Here and later on { 
$(u_0, w)_{\Sigma_0} := \int_{\Sigma_0} u_0(x) \, {w}(x,0) dx = \int_{\Omega} u_0(x) \,{w}(x,0) dx$.}
According to the well-establish arguments (see \cite{LMR:Ladyzhenskaya:1954, LMR:Ladyzhenskaya:1985,
LMR:Zeidler:1990a}), without loss of generality, we homogenise the problem, i.e., consider the problem 
with zero initial conditions $u_0 = 0$.

In order to be able to provide efficient discretisation method,
we introduce a stabilised weak formulation of \eqref{eq:equation} with time-upwind test functions
\begin{equation}
\lambda \, w + \mu \, \partial_t w, \quad 
w \in V^{\nabla_x \partial_t}_{0, \underline{0}} := \{ w \in V^{\Delta_x}_{0, \underline{0}} : \nabla_x \partial_t w \in \L{2}(Q) \}
, \quad 
\lambda, \mu \geq 0.
\label{eq:upwind-test}
\end{equation}
%
%
where $\lambda$ and $\mu$ are positive constants. We arrive at the following space-time 
formulation \cite{LMR:Hansbo:1994, LMR:Johnson:1987, LMR:JohnsonSaranen:1986}: find \linebreak
$u \in V_0$ satisfying 
\begin{equation}
a_s (u, w) = l_s (w), \quad \forall w \in V_{0},
\label{eq:stabilized-bilinear-form}
\end{equation}
where 
$$a_s (u, w) := \big(\partial_t u,  \lambda \, w + \mu \, \partial_t w \big)_Q 
+ \big(\nabla_x {u}, \nabla_x (\lambda \, w + \mu \, \partial_t w) \big)_Q$$
and 
$$l_s(w) := (f, \lambda \, w + \mu \, \partial_t w)_Q.$$

\begin{remark}
Notice that the approach presented in this paper (related to approximations and a posteriori error 
estimates) can be extended to more general parabolic equations, e.g., to those containing the term
$\dvrg_x(D(x,t) \nabla_x u(x,t))$ (where $D(x,t)$ is a positive definite matrix of diffusion coefficients) 
instead of $\laplace_x u(x,t)$ in (\ref{eq:equation}).
\end{remark}

Our main goal is to derive fully computable estimates for space IgA approximations of this class of the 
problems. For this purpose, we use the functional approach to a posteriori error estimates.
Initially, their simplest form have been obtained for a heat equation in \cite{LMR:Repin:2002}. In 
\cite{LMR:GaevskayaRepin:2005}, these estimates have been tested for generalised diffusion equation 
with a focus on an algorithmic part as well as on the comparison of two different forms of majorants. 
Evolutionary convection-diffusion equations and majorants for the approximations with jumps with 
respect to time variable has been considered in the paper \cite{LMR:RepinTomar:2010}. Finally, 
functional a 
posteriori error estimates for parabolic time-periodic BVPs as well as for optimal control problems have 
been studied in \cite{LMR:LangerRepinWolfmayr:2015} and \cite{LMR:LangerRepinWolfmayr:2016}.

This approach has been extended to a wide class of problems. Paper 
\cite{LMR:MatculevichRepin:2014} analyses the robustness of the majorant in the cases with drastic 
changes in values of the reaction parameter in reaction-diffusion time-dependent problems. It also 
discusses the quality of the indicator that follows from the functional error estimate. Moreover, it introduces 
the minorant of the error that is tested in several numerical examples and compared to the majorant. 
In order to make the estimates applicable to problems in domains with
complicated geometry, the domain decomposition technique in combination with local Poincar\'{e} 
inequalities can be used (see \cite{LMR:MatculevichNeitaanmakiRepin:2015}). 
The question on the constants in above-mentioned local inequalities is addressed in 
\cite{LMR:MatculevichRepin:2015}, where sharp bounds of them have been found for some simplexes
that are typically used in numerical methods. Numerical 
properties of above-mentioned error estimates w.r.t. the time-marching and space-time 
method are discussed in \cite{LMR:Matculevich:2015}.

\section{Overview of the IgA framework}
\label{eq:preliminaries}

For the convenience of the reader, we first recall the general concept of the IgA approach, the definition 
of B-splines (NURBS) and their use in the geometrical representation of the space-time cylinder $Q$,
as well as  
the construction of the IgA trial spaces, which are used to approximate solutions satisfying the variational 
formulation of \eqref{eq:variational-formulation}.

Let $p \geq 2$ be a polynomial degree and $n$ denote the number of basis functions used to construct 
a $B$-spline curve. The knot-vector in one dimension is a non-decreasing set of coordinates in the 
parameter domain, written as $\Xi = \{ \xi_1, \ldots, \xi_{n+p+1}\}$, $\xi_i \in \mathds{R}$, where 
$\xi_1 = 0$ and $\xi_{n+p+1} = 1$. The knots can be repeated, and the multiplicity of the $i$-th knot is 
indicated by $m_i$. Throughout the paper, we consider only open knot vectors, i.e., the multiplicity of the 
first and the last knots are equal to $m_1 = m_{n+p+1} = p+1$. 
For example, for the one-dimensional parametric domain $\Qhat = (0, 1)$, there is an underlying mesh 
${\mathcal{\Khat}}_h$ of elements $\Khat \in \mathcal{\Khat}_h$, such that each of them is 
constructed by the distinct neighbouring knots. The global size of $ \mathcal{\Khat}_h$ is denoted by 
$$\hhat := \max\limits_{\Khat \in \mathcal{\Khat}_h} \{ \hhat_{\Khat}\}, 
\quad \mbox{where} \quad 
\hhat_{\Khat} := {\rm diam} (\Khat).$$
For the time being, we assume locally quasi-uniform meshes, 
i.e., the ratio of two neighbouring elements $\Khat_i$ and $\Khat_j$ satisfies the inequality 
$c_1 \leq \tfrac{\hhat_{\Khat_i}}{\hhat_{\Khat_j}} \leq c_2$, where $c_1, c_2 > 0$.

The univariate B-spline basis functions $\Bhat_{i, p}: \Qhat \rightarrow \mathds{R}$ are defined by 
means of Cox-de Boor recursion formula as follows:
\begin{alignat}{2}
\Bhat_{i, p} (\xi) & := \tfrac{\xi - \xi_i}{\xi_{i+p} - \xi_i} \, \Bhat_{i, p-1} (\xi)
                         + \tfrac{\xi_{i+p+1} - \xi}{\xi_{i+p+1} - \xi_{i+1}} \Bhat_{i+1, p-1} (\xi), \quad 
\Bhat_{i, 0} (\xi) & \, := 
\begin{cases} 
1 & \mbox{if} \quad \xi_i \leq \xi \leq \xi_{i+1},  \\
0 & \mbox{otherwise},
\end{cases}
\label{eq:b-splines}
\end{alignat}
%
where a division by zero is defined to be zero. One of the most crucial properties of these basis functions 
\label{eq:b-splines} is that they are $(p-m_i)$-times continuously differentiable across the 
$i$-th knot with multiplicity $m_i$. Hence, if, for instance, $m_i = 1$ for every inner knot, B-splines of 
a degree $p$ are $C^{p-1}$ 
continuous. For the knots lying on the boundary of the parametric domain, where the multiplicity is 
$p + 1$, the B-spline is discontinuous ($C^{-1}$ function). We note, that analysis of this paper is
only concerned with a single-patch domain. The extension to the multi-patch case (the case, in which 
the physical domain is decomposed into several simple patches) will be a focus of a subsequent paper.

\begin{figure}
\centering
\begin{tikzpicture}[scale=0.6]
\def \ticksize {0.25};
\def \T {6.0};
\coordinate [label={above left:$a$}] (a) at (1.0, 0.0);
\coordinate [label={above right:$b$}] (b) at (5.0, 0.0);
\coordinate [label={above right:$$}] (aT) at (1.0, \T);
\coordinate [label={above right:$$}] (bT) at (5.0, \T);
\coordinate [label={below:$\Sigma_0$}] (Sigma_0) at (2.5, -0.2);
\coordinate [label={above:$\Sigma_T$}] (Sigma_T) at (3.5, \T + 0.2);
\coordinate [label={left:$T$}] (T) at (0.0, \T);
\coordinate [label={above left:$0$}] (O) at (0.0, 0.0);
\coordinate [label={above:$Q$}] (Q) at (3.0, \T/2);
\coordinate [label={above:$\Sigma$}] (sigma) at (5.3, 3.5);

\coordinate [label={above left:$$}] (xi10) at (2.0, 0.0);
\coordinate [label={above left:$$}] (xi20) at (3.5, 0.0);

\coordinate [label={above left:$$}] (xi1T) at (2.0, \T);
\coordinate [label={above left:$$}] (xi2T) at (3.5, \T);

\coordinate [label={above left:$$}] (yi10) at (1.0, 1.0);
\coordinate [label={above left:$$}] (yi20) at (5.0, 1.0);

\coordinate [label={above left:$$}] (yi12) at (1.0, 2.5);
\coordinate [label={above left:$$}] (yi22) at (5.0, 2.5);

\coordinate [label={above left:$$}] (yi13) at (1.0, 4.5);
\coordinate [label={above left:$$}] (yi23) at (5.0, 4.5);

\draw[->,thin,gray] (-1,0) --++(7,0)node[below left]{$x$};
\draw[->,thin,gray] (0,-1) --++(0,8)node[below left]{$t$};

\draw[thin, gray] (1, \ticksize) -- (1, -\ticksize);
\draw[thin, gray] (5.0, \ticksize) -- (5.0, -\ticksize);
\draw[thin, gray] (0.0-\ticksize, 6.0) -- (0.0 +\ticksize, 6.0);

\draw[black] (a) -- (b);
\draw[black] (aT) -- (bT);
\draw[black] (a) to [bend left=0] (aT);
\draw[black] (b) to [bend left=0] (bT);

\draw[black] (xi10) to [bend left=-0] (xi1T);
\draw[black] (xi20) to [bend left=0] (xi2T);

\draw[black] (yi10) -- (yi20);
\draw[black] (yi12) -- (yi22);
\draw[black] (yi13) -- (yi23);

\draw[->,thin,gray] (8,0) --(13,0)node[below left]{$\hat{x}$};
\draw[->,thin,gray] (9,-1) --(9,5)node[below left]{$\hat{t}$};

\draw[black] (9, 0) -- (12, 0) -- (12, 3) -- (9, 3) -- (9, 0);
\coordinate [label={above:$\Qhat$}] (Q) at (10.5, \T/4);

\coordinate [label={above left:$$}] (xi10hat) at (9.7, 0.0);
\coordinate [label={above left:$$}] (xi20hat) at (9.7, 3.0);

\coordinate [label={above left:$$}] (xi1That) at (11.1, 0.0);
\coordinate [label={above left:$$}] (xi2That) at (11.1, 3.0);

\coordinate [label={above left:$$}] (yi10hat) at (9.00, 0.6);
\coordinate [label={above left:$$}] (yi20hat) at (12.00, 0.6);

\coordinate [label={above left:$$}] (yi12hat) at (9.00, 1.5);
\coordinate [label={above left:$$}] (yi22hat) at (12.00, 1.5);

\coordinate [label={above left:$$}] (yi13hat) at (9.00, 2.5);
\coordinate [label={above left:$$}] (yi23hat) at (12.00, 2.5);

\draw[black] (xi10hat) -- (xi20hat);
\draw[black] (xi1That) -- (xi2That);

\draw[black] (yi10hat) -- (yi20hat);
\draw[black] (yi12hat) -- (yi22hat);
\draw[black] (yi13hat) -- (yi23hat);

\draw[black, ->] (5.5, \T/2) to [bend left=15] (8.0, \T/2);
\draw[black, <-] (5.5, \T/4) to [bend left=-15] (8.0, \T/4);

\draw (7, \T/2 + 0.3) node[label={above:{$\Phi^{-1}$}}] (Phiinverse) {};  
\draw (7, \T/4 - 0.3) node[label={below:{$\Phi$}}] (Phi) {};  

\coordinate [label={above left:$0$}] (Ohat) at (9.0, 0.0);
\coordinate [label={above right:$1$}] (1hat) at (12.0, 0.0);
\coordinate [label={above left:$1$}] (1hat) at (9.0, 3.0);

\draw[thin, gray] (9, \ticksize) -- (9, -\ticksize);
\draw[thin, gray] (12, \ticksize) -- (12, -\ticksize);
\draw[thin, gray] (-\ticksize + 9.0, 3.0) -- (9.0+\ticksize, 3.0);

\end{tikzpicture}
\caption{Mapping of the reference domain $\Qhat$ to non-moving in time cylinder $Q$. }
\end{figure}
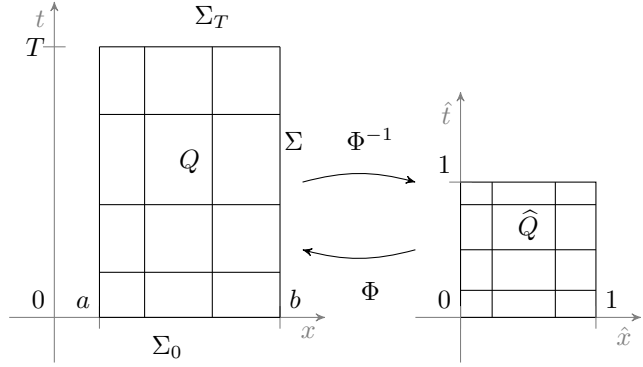

Consider also the multivariate B-splines on the space-time parameter domain $\Qhat := (0, 1)^{d+1}$, 
$d = \{1, 2, 3\}$, as a 
tensor-product of the corresponding univariate B-splines basis functions. 
For that, we define the knot vector dependent on the space-time direction 
$\Xi^\alpha = \{ \xi^\alpha_1, \ldots, \xi^\alpha_{n^\alpha+p^\alpha+1}\}$, 
$\xi^\alpha_i \in \mathds{R}$, where $\alpha = 1, \ldots, d+1$ is the index indicating the direction in 
space or time. 
Furthermore, we introduce 
the set of multi-indices 
$${\mathcal{I}} = \big\{\, i = (i_1,  \ldots, i_{d+1}): i_\alpha = 1, \ldots, n_\alpha; 
\alpha = 1, \ldots, d+1 \big\},$$
and multi-indices standing for the order of polynomials $p := (p_1, \ldots, p_{d+1})$. 
The tensor-product of the univariate B-spline basis functions generates multivariate B-spline basis 
functions
\begin{equation}
\Bhat_{i, p} ({\boldsymbol \xi}) 
:= \prod\limits_{\alpha = 1}^{d+1} \Bhat_{i_\alpha, p_\alpha} (\xi^\alpha), 
\quad \mbox{where} \quad {\boldsymbol \xi} = (\xi^1, \ldots, \xi^{d+1}) \in \Qhat.
\label{eq:mutlivariable}
\end{equation}
%
The univariate and multivariate NURBS basis functions are defined in the parametric domain by means 
of the corresponding B-spine basis functions $\big\{ \Bhat_{i, p} \big\}_{i \in {\mathcal{I}}}$. 
For given $p := (p_1, \ldots, p_{d+1})$ and for any $i \in {\mathcal{I}}$, we define the 
NURBS basis functions $\Rhat_{i, p}$ as follows:
\begin{equation}
\Rhat_{i, p}: \Qhat \rightarrow \mathds{R}, \quad 
\Rhat_{i, p} ({\boldsymbol \xi}) := \tfrac{w_i \, \Bhat_{i, p} ({\boldsymbol \xi})}{W({\boldsymbol \xi})},
\end{equation}
with a weighting function
\begin{equation}
W: \Qhat \rightarrow \mathds{R}, \quad 
W({\boldsymbol \xi}) := \sum_{i \in {\mathcal{I}}} w_i \, \Bhat_{i, p} ({\boldsymbol \xi}),
\end{equation}
where $w_i$ are positive real numbers. 

The physical space-time domain $Q \subset \mathds{R}^{d+1}$ is defined from the parametric domain 
$\Qhat = (0, 1)^{d+1}$ by the geometrical mapping:
\begin{equation}
\Phi: \Qhat \rightarrow Q := \Phi(\Qhat ) \subset \mathds{R}^{d+1}, \quad 
\Phi({\boldsymbol \xi}) := \sum_{i \in \mathcal{I}} \Rhat_{i, p}({\boldsymbol \xi}) \, {\bf P}_i,
\label{eq:geom-mapping}
\end{equation}
where 
$\{{\bf P}_i\}_{i \in \mathcal{I}} \in \mathds{R}^{d+1}$ are the control points. For simplicity, 
we assume below the same polynomial degree for all coordinate directions, i.e., 
$p_{\alpha} = p$ for all $\alpha = 1,  ... , d+1$.%

By means of latter geometrical mapping \eqref{eq:geom-mapping}, the physical mesh $\mathcal{K}_h$ 
is defined on the space-time domain $Q$, whose elements are images of elements of 
$\mathcal{\Khat}_h$, i.e., $\mathcal{K}_h := \big\{K = \Phi(\Khat) : \Khat \in \mathcal{\Khat}_h \big\}$.
The global mesh size is denoted by 
\begin{equation}
h := \max\limits_{K \in \mathcal{K}_h} \{ \, h_{K}\,\}, \quad  
h_{K} := \| \nabla \Phi \|_{\L{\infty} (K)} \hhat_{\Khat}.
\label{eq:global-mesh-size}
\end{equation}
Moreover, we assume that the physical mesh is also quasi-uniform, i.e., there exists a positive constant $C_u$ 
independent of $h$, such that
\begin{equation}
h_{K} \leq h \leq C_u \, h_{K}.
\label{eq:hk-and-h-relation}
\end{equation}
The discretisation spaces on $Q$ are constructed by a push-forward of the NURBS basis functions
\begin{equation}	
V_h := {\rm span} \,\big\{ \phi_{h, i} := \Rhat_{i, p} \circ \Phi^{-1} \big\}_{i \in \mathcal{I}},
\label{eq:vh-v0h}
\end{equation}
where we assume that the geometrical mapping $\Phi$ is invertible in $Q$, with smooth inverse on each 
element $K \in \mathcal{K}_h$ (see \cite{LMR:TagliabueDedeQuarteroni:2014} and 
\cite{LMR:Bazilevsetal2006} for more details). Moreover, we introduce the subspaces
$$
V_{0h} := V_h \cap V_{0, \underline{0}}
$$
for the functions fulfilling homogeneous boundary condition.
%

Let us recall two fundamental inequalities, i.e., scaled trace and inverse inequalities, that are 
important for the derivation of a priori error estimates of the space-time IgA scheme presented in the 
further sections.
\begin{lemma}{\rm \bf \cite[Theorem 3.2]{LMR:EvansHughes2013}}
\label{lm:lemma-2}
Let $K \in \mathcal{K}_h$, then the {\bf \em scaled trace inequality} 
\begin{equation}
\| v \|_{\partial K} \leq C_{tr} \, h_{K}^{-\rfrac{1}{2}} (\| v \|_{K} + h_{K} \, \| \nabla_x v \|_{K})
\end{equation}
holds for all $v \in \H{1}(K)$, where $h_{K}$ is the diameter of the element $K \in \mathcal{K}_h$, and 
$C_{tr}$ is a positive constant {independent of $K \in \mathcal{K}_h$}.
\end{lemma}

\begin{lemma}
{\rm \bf \cite[Theorem 4.1]{LMR:Bazilevsetal2006} and \cite[Theorem 4.2]{LMR:EvansHughes2013}}
\label{lm:lemma-3}
Let $K \in \mathcal{K}_h$, then the {\bf \em inverse inequalities}
\begin{equation}
\| \nabla v_h\|_{K} \leq C_{inv, 1} \,h_{K}^{-1} \, \| v_h\|_{K}
\label{eq:int-inequality-1} \\
\end{equation}
and 
\begin{equation}
\| v_h\|_{\partial K} \leq C_{inv, 0} \,h_{K}^{-\rfrac{1}{2}} \, \| v_h\|_{K}
\label{eq:int-inequality-0}
\end{equation}
%
%
hold for all $v_h \in V_{h}$, where $C_{inv, 0}$ and $C_{inv, 1}$ are positive  constants 
independent of $K \in \mathcal{K}_h$.

\begin{remark}
Inequality \eqref{eq:int-inequality-1} yields inequalities with partial derivatives in space and time
\begin{alignat}{2}
\| \partial_{x_i} v_h\|_{K} & \leq C_{inv, 1} \,h_{K}^{-1} \, \| v_h\|_{K}, \quad i = 1, \ldots, d, 
\quad \mbox{and} \quad 
\| \partial_t v_h\|_{K} & \leq C_{inv, 1} \,h_{K}^{-1} \, \| v_h\|_{K},
\end{alignat}
since $\nabla v_h = (\nabla_x v_h, \partial_t v_h)$. Nevertheless, for the anisotropic case (w.r.t. 
to {spatial} and time derivatives) the constant $C_{inv, 1}$ 
will depend on the direction, i.e., 
\begin{alignat}{2}
\| \partial_{x_i} v_h\|_{K} & \leq C^{x_i}_{inv, 1} \,h_{K}^{-1} \, \| v_h\|_{K} \quad \mbox{and} \quad 
\| \partial_t v_h\|_{K} & \leq C^{t}_{inv, 1} \,h_{K}^{-1} \, \| v_h\|_{K}.
\end{alignat}
\end{remark}

\end{lemma}

\begin{remark}
Due to the higher smoothness of basis functions ($p \geq 2$), the obtained approximations are 
generally $C^{p-1}$-continuous, provided that the inner knots have the multiplicity $1$. Moreover, basis 
functions of degree $p \geq 2$ are at least $C^1$-continuous if the multiplicity of the inner knots is less 
than or equal to $p-1$. This automatically implies that their gradients are in $H^{\dvrg_x, 0}(Q)$. 
Therefore, there is no need for constructing the projection of $\nabla_x u_h \in [ \L{2}(Q)]^d$ 
into $H^{\dvrg_x, 0}(Q)$. 
\end{remark}

\section{Discretisation of stabilized weak formulation}
\label{sec:discretization-non-moving-domain}

Stable space-time IgA scheme for parabolic equations have been presented and analysed in 
\cite{LMR:LangerMooreNeumueller:2016a}, where the authors proved its efficiency for fixed and 
moving spatial computational domains. In particular, it was shown that the corresponding discrete bilinear 
form is elliptic on the IgA space (w.r.t. a discrete energy norm), bounded, consistent. Moreover, 
the approximation results for the IgA spaces yields an a priori discretisation error estimate w.r.t. 
the same norm. In our work, we consider slightly modified energy norm and prove that the same 
properties of the scheme remain valid. 


We assume that the used spline-basis have a sufficiently high order, such that 
%
\begin{equation}
w_h \in V_{0h} \subset V_{0, \underline{0}}
%
\label{eq:h-partialtgradx-1}
\end{equation}
In order to derive a {\em stable discrete IgA space-time scheme} for (15),
we set the parameters in \eqref{eq:upwind-test} to $\lambda = 1$ and $\mu = \delta_{s, h} = \theta \, h$, 
where $\theta$ is a positive constant and $h$ is the mesh-size defined in \eqref{eq:global-mesh-size}, 
such that
\begin{equation}
w_h + \delta_{s, h} \, \partial_t w_h, \quad \delta_{s, h} = \theta h, \quad w_h \in V_{0h}. 
\label{eq:test-function}
\end{equation}
%
Hence, \eqref{eq:stabilized-bilinear-form} implies the discrete stabilised space-time problem: find
$u_h \in V_{0h}$ satisfying 
\begin{equation}
a_{s, h} (u_h, w_h) = l_h(w_h), \quad \forall w_h \in V_{0h},
\label{eq:discrete-scheme}
\end{equation}
where 
$$a_{s, h} (u_h, w_h) := (\partial_t u_h, w_h)_Q + \delta_{s, h} \,(\partial_t u_h, \partial_t w_h)_Q
+ (\nabla_x {u_h}, \nabla_x {w_h})_Q + \delta_{s, h} \ (\nabla_x {u_h}, \partial_t (\nabla_x w_h))_Q$$
and 
$$l_h(v_h) := (f, w_h + \delta_{s, h} \, \partial_w v_h)_Q.$$
$V_{0h}$-coercivity of $a_{s, h}(\cdot, \cdot): V_{0h} \times V_{0h} \rightarrow \mathds{R}$ 
w.r.t. the norm 
\begin{equation}
| \! |\! | v_h | \! |\! |^2_{s, h} 
:= \| \nabla_x {v_h }\|^2_{Q} 
+ \delta_{s, h} \, \| \partial_t v_h  \|^2_{Q} + 
\| v_h \|^2_{\Sigma_T} 
+ \delta_{s, h} \, \| \nabla_x v_h  \|^2_{\Sigma_T}, 
\label{eq:error-non-moving}
\end{equation}
follows from Lemma below.
\begin{lemma}
\label{lm:lemma-1}
The form $a_{s, h}(\cdot, \cdot): V_{0h} \times V_{0h} \rightarrow \mathds{R}$ is strongly $V_{0h}$-coercive w.r.t. the norm 
$| \! |\! | \cdot | \! |\! |^2_{s, h}$, i.e., there exists a positive constant $\mu_e$ such that
$$a_{s, h} (w_h, w_h) \geq \mu_e | \! |\! | w_h | \! |\! |^2_{s, h}, \quad \forall w_h \in V_{0h},$$
where $\mu_e = \tfrac{1}{2}$.
\end{lemma}
\ProofBegin
By considering $a_{s, h}(w, w)$, we arrive at 
$$a_{s, h} (w_h, w_h) = \lambda \, \| \nabla_x {w_h}\|^2_{Q} 
                          + \mu \, \| \partial_t w_h \|^2_{Q} 
                          + \tfrac{\lambda}{2} \, \| w_h \|^2_{\Sigma_T} 
                          + \tfrac{\mu}{2} \, \| \nabla_x w_h  \|^2_{\Sigma_T}
                  \geq \tfrac{1}{2} \, | \! |\! | w_h | \! |\! |^2_{s, h}.$$  
%
\ProofEnd

The latter property implies the existence and uniqueness of the discrete solution $u_h \in V_{0h}$, 
which can be written as
%
%
$$
u_h(x, t) = u_h(x_1, . . . , x_d, x_{d+1}) := \Sum_{i \in \mathcal{I}} \underline{u}_{h,i} \, \phi_{h,i}.
$$
Here, 
$\underline{u}_h := [ \underline{u}_{h,i}]_{i \in \mathcal{I}} \in {\mathds{R}}^{|\mathcal{I}|}$ is 
a vector of unknowns (degrees of freedom) defined by a system of equations
$$K_h \, \underline{u}_h = f_h$$
with the matrix $K_h := [K_{h,ij} = a_{s, h}(\phi_{h, i}, \phi_{h, j}) ]_{i, j \in \mathcal{I}}$ and the
right-hand side (RHS) $f_h := [f_{h,i} = l_h(\phi_{h, i}) ]_{i \in \mathcal{I}} \in {\mathds{R}}^{|\mathcal{I}|}$ (generated by IgA discretisation for elliptic problems). From Lemma \ref{lm:lemma-1}, it also follows that 
the matrix $K_h$ is regular (the condition number of $K_h$ is bounded by a constant independent 
of $h$).

\vskip 10pt
Several results, crucial for an a priori error estimation, can be shown using \cite[Lemma 2, 3]{LMR:LangerMooreNeumueller:2016a}. 
\begin{lemma}
\label{lm:lemma-4}
The bilinear form $a_{s, h}(\cdot, \cdot)$ is uniformly bounded on $V_{0h, *} \times V_{0h}$, 
where $V_{0h, *} := {{H^{2}_{0}}(Q) + V_{0h}}$, i.e., 
$$
|a_{s, h}(w, w_h)| \leq \mu_b \, | \! |\! | w | \! |\! |^2_{s, h, *} \, | \! |\! | w_h | \! |\! |^2_{s, h},  \quad 
\forall v \in V_{0h, *}, \forall v_h \in V_{0h}, $$
where  
$$| \! |\! | w | \! |\! |^2_{s, h, *} := | \! |\! | w | \! |\! |^2_h + \delta_{s, h}^{-1} \, \| w \|^2_{Q}, 
\quad \forall w \in V_{0h, *}, $$
and $\mu_b$ is a positive constant, independent of $h$.

\end{lemma}
\ProofBegin
The proof follows the lines of the proof in \cite{LMR:LangerMooreNeumueller:2016a}, 
where $|a_{s, h}(w, w_h)|$ is estimated term by term. First, we integrate $(\partial_t w, w_h)_Q$
by parts, i.e., 
%
\begin{equation}
(\partial_t w, w_h)_Q 
= (w, w_h)_{\Sigma_T} - (w, \partial_t  w_h)_Q,
\label{eq:int-by-part}
\end{equation}
and then estimate each of the terms in the RHS of \eqref{eq:int-by-part} by the H\"{o}lder 
inequality:
$$(w, w_h)_{\Sigma_T} - (v, \partial_t  v_h)_Q 
\leq \Big[\| w \|^2_{\Sigma_T}\Big]^{\rfrac{1}{2}} \, \Big[\| w_h \|^2_{\Sigma_T}\Big]^{\rfrac{1}{2}} 
+ \Big[\delta_{s, h}^{-1} \, \| w \|^2_Q\Big]^{\rfrac{1}{2}}\, 
\Big[\delta_{s, h} \, \|\partial_t w_h\|^2_Q\Big]^{\rfrac{1}{2}}.$$
%
The second and the third terms in $a_{s, h}(w, w_h)$ 
are estimated analogously, whereas for the forth term, 
\cite{LMR:LangerMooreNeumueller:2016a} suggests
applying the inverse inequality \eqref{eq:int-inequality-1} and 
the condition for quasi-uniform meshes introduced in \eqref{eq:hk-and-h-relation}. If we 
summarise all the estimates and use the fact that both $ \delta_{s, h} \, \| \nabla_x w_h\|_{\Sigma_T}$ and 
$\delta_{s, h} \, \| \nabla_x w\|_{\Sigma_T}$ are positive quantities (which can be added to the RHS), 
we arrive at the relation
\begin{alignat*}{2}
|a_{s, h}(w, w_h)| & \leq 
\Big[\| w \|^2_{\Sigma_T} 
+ \delta_{s, h}^{-1} \, \|w\|^2_Q 
+ \delta_{s, h} \, \| \partial_t w \|^2_Q 
+ 2\, \| \nabla_x {w}\|^2_Q\Big]^{\rfrac{1}{2}} \\
& \qquad \qquad 
\Big[\| w_h \|^2_{\Sigma_T} 
+ \delta_{s, h} \, \|\partial_t w_h\|^2_Q 
+ \delta_{s, h} \, \| \partial_t w_h \|^2_Q 
+ \| \nabla_x {w_h}\|^2_Q 
+ C_u^2 \, C^2_{inv, 1} \, \theta^2 \| \nabla_x w_h \|^2_Q \Big]^{\rfrac{1}{2}} \\
& \leq 
\Big[\| w \|^2_{\Sigma_T} + \delta_{s, h} \, \| \nabla_x w \|^2_{\Sigma_T}
+ \delta_{s, h}^{-1} \, \|w\|^2_Q 
+ \delta_{s, h} \, \| \partial_t w \|^2_Q 
+ 2\, \| \nabla_x {w}\|^2_Q \Big]^{\rfrac{1}{2}} \, \\
& \qquad \qquad
\Big[\| w_h \|^2_{\Sigma_T} + \delta_{s, h} \, \| \nabla_x {w_h} \|^2_{\Sigma_T}
+ 2 \, \delta_{s, h}\, \|\partial_t w_h\|^2_Q 
+ (1 + C_u^2 \, C^2_{inv, 1} \, \theta^2) \, \| \nabla_x {w_h}\|^2_Q \Big]^{\rfrac{1}{2}} \\
& \leq \mu_b \, | \! |\! | w  | \! |\! |_{s, h, *} \, | \! |\! | w_h | \! |\! |_{h},
\end{alignat*}
where 
$\mu_b 
= \Big(\max \Big\{ 1 + C_u^2 \, C^2_{inv, 1} \, \theta^2, 2 \Big\}\Big)^{\rfrac{1}{2}}.$ 
\ProofEnd
\vskip 10pt

For the completeness, it is worth recalling basic results on the approximation properties of spaces 
generated by B-splines (NURBS) that follow from 
\cite[Section 3]{LMR:Bazilevsetal2006} and \cite[Section 4]{LMR:TagliabueDedeQuarteroni:2014}. They 
state the existences of a projection operator $\Pi_h: H^s_{0, \underline{0}}(Q) \rightarrow V_{h}$, 
$s \in \mathds{N}, s \geq 0$, that provide the 
asymptotically optimal approximation result.
\begin{lemma}
\label{lm:lemma-5}
Let $l, s \in \mathds{N}$ be $0 \leq l \leq s \leq p+1$, and $w \in H^s_{0, \underline{0}}(Q)$. 
Then, there exists a projection operator \linebreak 
$\Pi_h: H^s_{0, \underline{0}}(Q) \rightarrow V_{0h}$ and a positive constant $C_s$ such 
that 
\begin{equation}
\sum_{K \in \mathcal{K}_h} | w - \Pi_h w|^2_{H^l(K)}
 \leq C^2_s h^{2(s-l)} \| w \|^2_{\H{s}(Q)},
 \label{eq:c-s-interpolation-estimate}
 \end{equation}
 where $h$ is a global mesh size defined by (\ref{eq:global-mesh-size}) and $C_s$ is a constant 
 only dependent on degrees $s, l,$ and $p$, the shape regularity of $Q$, described by $\Phi$ and 
 its gradient. 
 \end{lemma}
 \ProofBegin
 The proof follows the lines of \cite[Subsection 3.3, 3.4]{LMR:Bazilevsetal2006}, 
 \cite[Proposition 3.1]{LMR:TagliabueDedeQuarteroni:2014}, and 
 \cite[Chapter 4]{LMR:BeiraodaVeigaBuffaRivasSangalli:2011}.
\ProofEnd
\vskip 10pt

If the multiplicity of each inner knot $m_i \leq p + 1 - l$, $\Pi_h w$ belongs to $\HD{s}{0, \underline{0}}(Q)$.
Then, Lemma \ref{lm:lemma-5} yields the global estimate
\begin{equation}
| w - \Pi_h w|_{H^l(Q)} \leq C_s h^{(s-l)} \| w \|_{\H{s}(Q)}.
\label{eq:global-approximtion-estimate}
\end{equation}
Both interpolation estimates 
\eqref{eq:c-s-interpolation-estimate} and \eqref{eq:global-approximtion-estimate} yield a priori 
estimates of the interpolation error $e_h =  w - \Pi_h w$, measured in terms of the $\L{2}$-norm and
the discrete norms $| \! |\! | \cdot  | \! |\! |_{h}$ and $| \! |\! | \cdot  | \! |\! |_{h, *}$, which 
we later need in order to obtain an a priori estimate for the discretisation error $u - u_h$.

\begin{lemma}
\label{lm:lemma-6}
Let $l, s \in \mathds{N}$ be $1 \leq l \leq s \leq p+1$, and $v \in \HD{s}{0, \underline{0}}(Q)$. Then, 
there exists a projection operator \linebreak $\Pi_h: \HD{s}{0, \underline{0}}(Q) \rightarrow V_{0h}$ (see Lemma \ref{lm:lemma-5}) 
and positive constants $C_1, C_2, C_3, C_4$, such that the following a priori error estimates hold
 \begin{alignat}{2}
 \| w - \Pi_h w \|_{\partial Q} & \leq C_1 \,h^{s - \rfrac{1}{2}} \| w \|_{\H{s}(Q)}, \label{eq:int-1}\\
 \delta^{\rfrac{1}{2}}_h 
 \| \nabla_x(w - \Pi_h w) \|_{\partial Q} & 
 \leq C_2 \,h^{s - 1} \| w \|_{\H{s}(Q)}, \label{eq:int-2}\\ 
 | \! |\! | w - \Pi_h w  | \! |\! |_{h} & \leq C_3 \,h^{s - 1} \| w \|_{\H{s}(Q)}, \label{eq:int-3}\\
| \! |\! | w - \Pi_h w  | \! |\! |_{s, h, *} & \leq C_4 \,h^{s - 1} \| w \|_{\H{s}(Q)}. \label{eq:int-4}
 \end{alignat}
\end{lemma}
\ProofBegin
Estimate \eqref{eq:int-1} follows straightforwardly from the proof of Lemma 6 in \cite{LMR:LangerMooreNeumueller:2016a}. Let us show that estimate \eqref{eq:int-2} holds:
\begin{alignat*}{2}
\delta_{s, h} \| \nabla_x (w - \Pi_h w) \|^2_{\partial Q}
& \leq \theta \, h  \, \Sum_{K \in \mathcal{K}_h} \Sum_{i = 1, ..., d} \, 
      \| \partial_{x_i} (w - \Pi_h w) \|^2_{\partial K \cap \partial Q} \\
& \mnote{\eqref{eq:int-inequality-0}}
\leq 2 \, \theta \, h \, \Sum_{K \in \mathcal{K}_h} \Sum_{i = 1, ..., d} \, 
       \Big( C^2_{0, inv} \, h^{-1}_{K} \, \Big(\| \partial_{x_i} (w - \Pi_h w) \|^2_{K} 
                                             + h^2_{K} \, | \partial_{x_i} (w - \Pi_h w) |^2_{K}\Big) \Big) \\
& \mnote{\eqref{eq:hk-and-h-relation}} \leq 2 \, \theta \, h \, \Sum_{K \in \mathcal{K}_h} \Sum_{i = 1, ..., d} \, 
       \Big( C^2_{0, inv} \, C_u \, h^{-1}\, \Big(\| \partial_{x_i} (w - \Pi_h w) \|^2_{K} 
                                             + h^2 \, \Sum_{j = 1, \ldots, d+1}| \partial_{x_i} \partial_{x_j} (w - \Pi_h w) |^2_{{K}}\Big) \Big) \\
& \leq 2 \, C^2_{0, inv} \, C_u \, \theta \,
         \Big( \|\nabla_{x} (w - \Pi_h w) \|^2_{Q} 
                                             + h^2 \Sum_{j = 1, \ldots, d+1}| \partial_{x_i} \partial_{x_j} (w - \Pi_h w) |^2_{{Q}}\Big) \\
& \mnote{\eqref{eq:c-s-interpolation-estimate}} 
   \leq 2 \, C_u \, \theta \, (1 + \, d\, (d+1)) \,  C^2_{0, inv} \, C^2_s \, h^{2(s-1)} \, \| w \|^2_{\H{s}(Q)}.
\end{alignat*}
Since $\| \nabla_x (w - \Pi_h w) \|^2_{\Sigma_T} \leq \| \nabla_x (w - \Pi_h w) \|^2_{\partial Q},$
the relations \eqref{eq:int-3} and \eqref{eq:int-4} hold.
\ProofEnd


\begin{lemma}
\label{lm:lemma-7}
If the solution $u \in \HD{1, 0}{0}(Q) \cap \HD{2}{} (Q)$, then it satisfies the consistency identity $a_{s, h}(u, v_h) = l_h(w_h)$,  
$\forall w_h \in V_{0h}$. 
\end{lemma}
\ProofBegin
The proof follows the steps of the lines of Lemma 7 in \cite{LMR:LangerMooreNeumueller:2016a}.
\ProofEnd

\vskip 10pt
\noindent
The main result of this section, namely, the a priori error estimate in the discrete norm $\| \cdot \|_h$, 
is formulated in the following theorem. 

\begin{theorem}
\label{lm:theorem-8}
Let $u \in H^{s}_{0}(Q) := H^{s}(Q) \cap H^{1, 0}_{0}(Q)$ with $s \in \mathds{N}$, $s \geq 2$, 
be the exact solution of \eqref{eq:variational-formulation}, and 
let $u_h \in V_{0h}$ be a solution of the space-time IgA scheme \eqref{eq:discrete-scheme} with some 
fixed parameter $\theta$. Then, the discretisation error estimate
\begin{equation}
\| u - u_h\|_h \leq C \, h^{r-1}\, \| u \|_{H^r(Q)} 
\label{eq:error-in-h}
\end{equation}
holds, where $C$ is a constant independent of $h$ and $r = \min \{ s, p+1 \}$.
\end{theorem}

\section{Error majorant}
\label{eq:general-error-estimate}

In this section, we derive error majorants of the functional type for stabilised weak formulations of 
parabolic I-BVPs with time-upwind test functions. 
These error estimates are used to obtain a posteriori error estimates for the distance 
$e = u - v$ between $u \in V^{\Delta_x}_{0, \underline{0}}$ and any 
$v \in V^{\Delta_x}_{0, \underline{0}} (V^{\nabla_x \partial_t}_{0, \underline{0}})$
(in particular, approximations produced by the space-time IgA method presented in the previous 
section) measured in terms of the norm
\begin{equation}
| \! |\! | e | \! |\! |^2_{s, \nu_i}
:= \nu_{1} \, \| \nabla_x {e}\|^2_{Q} +  \nu_{2} \, \| \partial_t e \|^2_{Q} 
+ \nu_{3} \, \| \nabla_x e \|^2_{\Sigma_T} + \nu_{4} \, \| e \|^2_{\Sigma_T},
\label{eq:error-general}
\end{equation}
where $\nu_{i} >0$, $i = 1, \ldots, 4$, are some weights (introduced throughout the derivation process). 

To obtain guaranteed error bounds of \eqref{eq:error-general}, we 
apply a method similar to the one 
developed in \cite{LMR:Repin:2002,LMR:MatculevichRepin:2014} for parabolic I-BVPs. 
For the derivation
process, we consider space of smoother functions $V^{\nabla_x \partial_t}_{0, \underline{0}}$ 
(cf. \eqref{eq:upwind-test}) equipped with the 
norm 
$$\| w \|_{V^{\nabla_x \partial_t}_{0, \underline{0}}} 
:= \sup\limits_{t \in [0, T]} \| \nabla_x w(\cdot, t)\|^2_Q + \| w \|^2_{V^{\Delta_x}_{0, \underline{0}}},$$
where 
$$\| w \|^2_{V^{\Delta_x}_{0, \underline{0}}} := \| \Delta_x w \|^2_Q + \| \partial_t w \|^2_Q,$$
which is dense in $V^{\Delta_x}_{0, \underline{0}}$. According to 
\cite[Remark 2.2]{LMR:Ladyzhenskaya:1985}, norms
$\| \cdot \|_{V^{\nabla_x \partial_t}_{0, \underline{0}}} \approx \| \cdot \|_{V^{\Delta_x}_{0, \underline{0}}}$.

Consider the sequence 
$u_n \in V^{\nabla_x \partial_t}_{0, \underline{0}}$.Then, the corresponding stabilised identity is formulated
as follows:
\begin{equation}
a_s (u_n, w) = (f_n, \lambda \, w + \mu \, \partial_t w)_Q, \quad \mbox{where} \quad
f_n = {(u_n)}_t - \Delta_x u_n \in \L{2}(Q).
\label{eq:stabilized-bilinear-form-short}
\end{equation}
%
%
{By subtracting $a_s(v_n, w)$, $v_n \in V^{\nabla_x \partial_t}_{0, \underline{0}}$,  
from \eqref{eq:stabilized-bilinear-form-short},
and by setting $w = e_n = u_n - v_n \in V^{\nabla_x \partial_t}_{0, \underline{0}}$, 
we arrive at the so-called `error-identity'}
%
%
%
\begin{multline*}
	\lambda \, \| \nabla_x {e_n} \|^2_Q  + \mu \, \| \, \partial_t e_n \|^2_Q  
	+ \tfrac12 \, (\mu \, \| \nabla_x{e_n} \|^2_{\Sigma_T} + \lambda \| e_n \|^2_{\Sigma_T}) \\
	= \lambda \Big( (f_n - \partial_t v_n,  e_n)_Q - (\nabla_x {v_n}, \nabla_x {e_n})_{Q} \Big) 
	+ \mu \Big( (f_n - \partial_t  v_n, \partial_t  e_n)_Q - (\nabla_x {v_n}, \nabla_x \, \partial_t e_n)_{Q} \Big),
\end{multline*}
%
which is used in the derivation of the majorants of \eqref{eq:error-general}
in Theorems \ref{th:theorem-majorant-general-1} and \ref{th:theorem-majorant-general-2}.
%
%
\begin{theorem}
\label{th:theorem-majorant-general-1}
%
For any $v \in V^{\Delta_x}_0$ and ${\flux} \in H^{\dvrg_x, 0}(Q)$, 
the following estimate holds:
%
\begin{alignat}{2}
\!\!\!
(2- \tfrac{1}{\gamma}) & 
  (\lambda\,\| \nabla_x e \|^2_{Q} + \mu  \, \| \partial_t e \|^2_{Q}) 
+ \mu \, \| \nabla_x{e} \|^2_{\Sigma_T} + \lambda \| e \|^2_{\Sigma_T} 
\leq \overline{\rm M}^{\rm I} (v, \flux; \gamma, \alpha_i) \nonumber \\[5pt]
& := \gamma\Big\{ \lambda \Big((1 + \alpha_1)\,\| \R_{\rm d} \|_{Q}^2 
      + (1 + \tfrac{1}{\alpha_1})\, \CFriedrichs^2 \, \|\R_{\rm eq}\|^2_{Q} \Big)
			+ \mu \Big( 
			(1 + \alpha_2) \, \| \dvrg_x \R_{\rm d}\|^2_{Q} 
			+ (1 + \tfrac{1}{\alpha_2})\,\| \R_{\rm eq} \|^2_{Q}
			\Big) \Big\},
\label{eq:estimate-1-non-moving}
\end{alignat}
where $\R_{\rm eq}$ and $\R_{\rm d}$ are defined by relations 
\begin{alignat}{2}
	\R_{\rm eq}  (v, \flux) & := f - \partial_t v + \dvrg_x \: \flux \quad \mbox{and} \quad 
	\R_{\rm d}  (v, \flux) := \flux - \nabla_x {v},
\label{eq:residuals}
\end{alignat}
$\CFriedrichs$ is the Friedrichs constant, 
$\lambda$ and $\mu$ are positive weights introduced in \eqref{eq:upwind-test},  
$\gamma \in \big[\tfrac{1}{2}, +\infty)$, and $\alpha_i > 0, i = 1, 2.$ 
\end{theorem}
%
%

\noindent{\bf Proof:} The RHS of the error-identity is modified by means of the relation 
$$(\dvrg_x \flux, \lambda \, e_n  \!+\! \mu \, \partial_t e_n)_Q 
+ \big(\flux, \nabla_x (\lambda \, e_n  \!+\! \mu \, \partial_t e_n)\big)_Q = 0.$$
The obtained result can be presented as follows:
\begin{multline}
	\lambda \, \| \nabla_x {e_n} \|^2_Q 
	+ \mu \, \| \partial_t \, e_n \|^2_Q
	+ \tfrac12 \, ( \mu \, \| \nabla_x{e_n} \|^2_{\Sigma_T} + \lambda \, \| e_n \|^2_{\Sigma_T}) \\
	= \lambda \, \big(\left(f_n - \partial_t v_n + \dvrg_x \: \flux,  e_n \right)_Q + ( \flux - \nabla_x {v_n}, \nabla_x {e_n})_{Q} \big) \\
	      + \mu \, \big( \left(f_n - \partial_t v_n + \dvrg_x \: \flux, \partial_t {e_n} \right)_Q + ( \flux - \nabla_x {v}, \nabla_x \partial_t {e_n})_{Q}\big).
	\label{eq:energy-balance-equation-with-flux}
\end{multline}
%
We proceed further by integrating by parts the term $ (\R_{\rm d}, \nabla_x \partial_t {e_n})_{Q}$:
\begin{equation*}
\mu\, \big(\R_{\rm d}, \nabla_x (\partial_t e_n)\big)_Q 
= \mu\, (\R_{\rm d}, \vectorn_x \, \partial_t e_n)_\Sigma
- \mu \, (\dvrg_x ( \flux - \nabla_x {v_n}), \partial_t e_n)_{Q}
= - \mu \, (\dvrg_x \flux - \Delta_x v_n, \partial_t e_n)_{Q}.
\end{equation*}
%
%
Using density arguments, i.e., 
$u_n \rightarrow u$,  
$v_n \rightarrow v \in V^{\Delta_x}_{0, \underline{0}}$, and 
$f_n \rightarrow f \in \L{2}(Q)$ for $n \rightarrow \infty$, 
we arrive at the identity formulated for $e = u - v$ with $u, v \in V^{\Delta_x}_{0, \underline{0}}$, i.e., 
\begin{multline}
	\lambda \, \| \nabla_x {e} \|^2_Q 
	+ \mu \, \| \partial_t \, e \|^2_Q
	+ \tfrac12 \, ( \mu \, \| \nabla_x{e} \|^2_{\Sigma_T} + \lambda \, \| e \|^2_{\Sigma_T}) \\
	= \lambda \, \big(\left(\R_{\rm eq},  e \right)_Q + (\R_{\rm d}, \nabla_x {e})_{Q} \big)
	      + \mu \, \big( \left(\R_{\rm eq}, \partial_t {e} \right)_Q - \mu \, (\dvrg_x \R_{\rm d}, \partial_t e)_{Q} \big).
	\label{eq:energy-balance-equation-with-flux-new}
\end{multline}
The first term on the RHS of \eqref{eq:energy-balance-equation-with-flux-new} is estimated by the 
H\"{o}lder and Friedrichs inequalities
\begin{alignat*}{2}
\lambda \, (\R_{\rm eq}, \, e)_Q         & \leq \CFriedrichs \, \lambda \, \|\R_{\rm eq}\|_{Q} \, \| \nabla_x e \|_{Q}. 
\end{alignat*}
The second, third, and forth terms can be treated analogously. 
Therefore, \eqref{eq:energy-balance-equation-with-flux-new} yields the estimate
\begin{alignat}{2}
\lambda \, \| \nabla_x{e} \|^2_Q 
\!+ \!\mu \, \| \partial_t \, e \|^2_Q
& + \tfrac12 \, ( \mu \, \| \nabla_x{e} \|^2_{\Sigma_T} \! 
   + \!\lambda \| e \|^2_{\Sigma_T}) 
\nonumber\\[5pt]
& \leq \, \lambda\, (\| \R_{\rm d} \|_{Q} \!+\! \CFriedrichs \, \|\R_{\rm eq}\|_{Q})\, 
  \| \nabla_x e \|_{Q}
\!+\! \mu \, 
  \big(\| \dvrg_x \R_{\rm d}\|_{Q} + \| \R_{\rm eq} \|_{Q}\big) \, 
        \| \partial_t e \|_{Q} \nonumber\\[5pt]
& \leq \, \big( \lambda \, (\| \R_{\rm d} \|_{Q} + \CFriedrichs \, \|\R_{\rm eq}\|_{Q})^2 + 
  \mu \,(\| \dvrg_x \R_{\rm d}\|_{Q} + \| \R_{\rm eq} \|_{Q})^2 \big)^{\rfrac{1}{2}}
  (\lambda \, \| \nabla_x e \|^2_{Q} + \mu \, \| \partial_t e \|^2_{Q})^{\rfrac{1}{2}}.
\label{eq:energy-balance-equation-3-general-estimate}
\end{alignat}
In order to regroup the terms on the RHS of \eqref{eq:energy-balance-equation-3-general-estimate}, 
we apply the Young inequality with positive scalar-valued parameters
$\gamma$, $\alpha_1$, and $\alpha_2$ and deduce the estimate
\begin{alignat*}{2}
\big( \lambda\, ( \| \R_{\rm d} \|_{Q} & + \CFriedrichs \, \|\R_{\rm eq}\|_{Q})^2 
  + \mu \,(\| \dvrg_x \R_{\rm d}\|_{Q} + \| \R_{\rm eq} \|_{Q})^2 \big)^{\rfrac{1}{2}}
  (\lambda \, \| \nabla_x e \|^2_{Q} + \mu \, \| \partial_t e \|^2_{Q})^{\rfrac{1}{2}} \\[5pt]
& \leq \tfrac{\gamma}{2} \, 
\big( \lambda\, (\| \R_{\rm d} \|_{Q} + \CFriedrichs \, \| \R_{\rm eq}\|_{Q})^2
      + \mu\, (\| \dvrg_x \R_{\rm d}\|_{Q} + \| \R_{\rm eq} \|_{Q})^2 \big)
+ \tfrac{1}{2\,\gamma} \, ( \lambda\, \| \nabla_x e \|^2_{Q} + \mu \, \| \partial_t e \|^2_{Q}) \\[5pt]
& \leq \tfrac{\gamma}{2} \, 
 \Big( \lambda \, \big( (1 + \alpha_1)\,\| \R_{\rm d} \|_{Q}^2 + (1 + \tfrac{1}{\alpha_1})\, \CFriedrichs^2 \, \|\R_{\rm eq}\|^2_{Q}\Big)
+ \mu\, \big( (1 + \alpha_2)\, \| \dvrg_x \R_{\rm d}\|^2_{Q} + (1 + \tfrac{1}{\alpha_2})\,\| \R_{\rm eq} \|^2_{Q}\big) \Big) \\[5pt]
& \qquad \qquad \qquad \qquad \qquad \qquad \qquad \qquad \qquad \qquad \qquad \qquad \qquad \qquad
+ \tfrac{1}{2\,\gamma} \, (\lambda\, \| \nabla_x e \|^2_{Q} + \mu  \, \| \partial_t e \|^2_{Q}),
\end{alignat*} 
where $\gamma$, $\alpha_1$, and $\alpha_2$ are positive scalar-valued parameters. Then, 
the obtained inequality yields \eqref{eq:estimate-1-non-moving}.
\hfill $\square$
\vskip 10pt

The next Theorem require higher regularity for both $v$ and ${\flux}$ with respect to time.
\begin{theorem}
\label{th:theorem-majorant-general-2}
For any $v \in V^{\nabla_x \partial_t}_{0, \underline{0}}$ and any ${\flux} \in H^{\dvrg_x, 1}(Q)$, 
the following inequality holds:
\begin{alignat}{2}
\!\!\!\!
(2 - \tfrac{1}{\zeta}&) (\lambda \,\| \nabla_x {e}\|^2_{Q} + \mu \, \| \partial_t e \|^2_{Q}) 
+ \mu \, (1 - \tfrac{1}{\epsilon}) \| \nabla_x e \|^2_{\Sigma_T} 
+ \lambda\, \| e \|^2_{\Sigma_T} \leq \overline{\rm M}^{\rm II}(v, \flux; \zeta, \beta_i, \epsilon) 
:= \epsilon \, \mu \| \R_{\rm d} \|^2_{\Sigma_T} \nonumber\\[5pt]
& \qquad \quad
+ \zeta
  \Big( \lambda \big( (1 + \beta_1)\big(
                                         (1 + \beta_2) \, \| \R_{\rm d} \|^2_{Q} 
				+ (1 + \tfrac{1}{\beta_2}) \CFriedrichs^2 \, \|\R_{\rm eq}\|^2_{Q} \big) 
				+ (1 + \tfrac{1}{\beta_1}) \, \tfrac{\mu^2}{\lambda^2} \, 
				   \| \partial_t \R_{\rm d} \|^2_Q\big)
         + \mu \, \| \R_{\rm eq} \|_{Q}^2 \Big),
\label{eq:estimate-2-non-moving-domain}
\end{alignat}
where $\R_{\rm eq} (v, y)$ and
$\R_{\rm d}  (v, y)$ are defined in \eqref{eq:residuals}, 
$\CFriedrichs$ is the Friedrichs constant, 
$\lambda$ and $\mu$ are positive weights introduced in 
\eqref{eq:upwind-test}, $\zeta \in \big[\tfrac{1}{2}, +\infty)$, $\epsilon \in [1, +\infty)$, 
and $\beta_i >0, i = 1, 2$.
\end{theorem}
%
\ProofBegin 
Now, we use a different transformation of the last term in the RHS of
\eqref{eq:energy-balance-equation-with-flux}:
\begin{equation}
\mu \, \big (\R_{\rm d}, \nabla_x (\partial_t e_n)\big)_Q
= \mu \, (\R_{\rm d}, \nabla_x e_n \, n_t)_{\Sigma_T}
- \mu \, (\partial_t \R_{\rm d}, \nabla_x e_n)_Q, 
\label{eq:last-term-general-estimate-2}
\end{equation}
where $n_t\big|_{\Sigma_T} = 1$. Analogously to the proof of 
Theorem \ref{th:theorem-majorant-general-1}, we use density arguments to obtain
\begin{multline}
	\lambda \, \| \nabla_x {e} \|^2_Q 
	+ \mu \, \| \partial_t \, e \|^2_Q
	+ \tfrac12 \, ( \mu \, \| \nabla_x{e} \|^2_{\Sigma_T} + \lambda \, \| e \|^2_{\Sigma_T}) \\
	= \lambda \, \big(\left(\R_{\rm eq},  e \right)_Q + (\R_{\rm d}, \nabla_x {e})_{Q} \big)
	      + \mu \, \big( (\R_{\rm eq}, \partial_t {e} )_Q 
	       \, (\R_{\rm d}, \nabla_x e \, n_t)_{\Sigma_T}
                - \, (\partial_t \R_{\rm d}, \nabla_x e)_Q \big).
	\label{eq:energy-balance-equation-with-flux-new}
\end{multline}
Since
%
\begin{alignat*}{2}
\mu \, (\R_{\rm d}, \nabla_x e)_{\Sigma_T} &
\leq  \, \tfrac{\mu}{2} \, ( \tfrac{1}{\epsilon} \, \| \nabla_x e \|^2_{\Sigma_T} 
     + \epsilon \| \R_{\rm d} \|^2_{\Sigma_T}), 
\quad \epsilon >0,
\end{alignat*}
and 
\begin{alignat*}{2}
- \mu \, (\partial_t \R_{\rm d}, \nabla_x e)_Q & 
\leq \, \mu \, \| \partial_t \R_{\rm d} \|_{Q} \, \| \nabla_x e \|_{Q},
\end{alignat*}
%
%
we obtain 
%
\begin{alignat*}{2}
\lambda \, & \| \nabla_x {e}\|^2_{Q} + \mu \, \| \partial_t e \|^2_{Q}
+ \tfrac{1}{2} \, ( \lambda \, \| e \|^2_{\Sigma_T} + \mu \, \| \nabla_x e \|^2_{\Sigma_T}) 
\nonumber\\[5pt]
& \leq \tfrac{\mu}{2}\,
\big(\tfrac{1}{\, \epsilon} \| \nabla_x e \|^2_{\Sigma_T} 
 + \epsilon \| \R_{\rm d} \|^2_{\Sigma_T}\big)
 + \big (\lambda\,  (\| \R_{\rm d} \|_{Q} 
                           + \CFriedrichs \, \|\R_{\rm eq}\|_{Q}) 
           + \mu \, \| \partial_t \R_{\rm d} \|_Q \big) \, \| \nabla_x e \|_{Q} 
+ \mu \, \| \R_{\rm eq} \|_{Q}\, \| \partial_t e \|_{Q} 
\nonumber\\[5pt]
& \leq \tfrac{\mu}{2} \,
\big(\tfrac{1}{\, \epsilon} \| \nabla_x e \|^2_{\Sigma_T} 
 + \epsilon \| \R_{\rm d} \|^2_{\Sigma_T}\big) 
 + \Big( \lambda \, \big( \| \R_{\rm d} \|_{Q} 
                           + \CFriedrichs \, \|\R_{\rm eq}\|_{Q}
                           + \tfrac{\mu}{\lambda} \, \| \partial_t \R_{\rm d} \|_Q
            \big)^2 
           + \mu \, \| \R_{\rm eq} \|_{Q}^2 
     \Big)^{\rfrac{1}{2}}
     (\lambda \, \| \nabla_x e \|^2_{Q} 
     + \mu \, \| \partial_t e \|^2_{Q})^{\rfrac{1}{2}} \nonumber\\[5pt]
& \leq \tfrac{\mu}{2} \,
\big(\tfrac{1}{\, \epsilon} \| \nabla_x e \|^2_{\Sigma_T} 
+ \epsilon \| \R_{\rm d} \|^2_{\Sigma_T}\big) 
+ \tfrac{1}{2\,\zeta}	
(\lambda \, \| \nabla_x e \|^2_{Q} + \mu \, \| \partial_t e \|^2_{Q})
\nonumber\\[5pt]
&  \qquad \qquad \qquad \qquad 
+ \tfrac{\zeta}{2} \, 
  \bigg(  \lambda \, 
           \Big( (1 + \beta_1) \,  \,
                 \big( (1 + \beta_2) \, \| \R_{\rm d} \|^2_{Q} 
	     + (1 + \tfrac{1}{\beta_2}) \CFriedrichs^2 \, \|\R_{\rm eq}\|^2_{Q} \big) 
	 + (1 + \tfrac{1}{\beta_1})  \, \tfrac{\mu^2}{\lambda^2} \, \| \partial_t \R_{\rm d} \|^2_Q \Big)
      + \mu \, \| \R_{\rm eq} \|_{Q}^2 \bigg). 
\end{alignat*}
%
This estimate yields \eqref{eq:estimate-2-non-moving-domain}.
%
\ProofEnd 
\vskip 10pt

\noindent
In Corollaries below, we consider a particular case related to the choice $\lambda = 1$ and 
$\mu = \delta_{s, h}$.
\begin{corollary}
\label{cor:majorant-1}
Assume that $v \in V^{\Delta_x}_{0, \underline{0}}$ 
%
and ${\flux} \in H^{\dvrg_x, 0}(Q)$. Then, Theorems \ref{th:theorem-majorant-general-1} yields 
the estimate
\begin{alignat}{2}
(2 - \tfrac{1}{\gamma}) & (\,\| \nabla_x {e}\|^2_{Q} + \delta_{s, h} \, \| \partial_t e \|^2_{Q}) 
+ \delta_{s, h} \, \| \nabla_x e \|^2_{\Sigma_T} 
+ \| e \|^2_{\Sigma_T} 
\leq \overline{\rm M}^{\rm I}_{\delta_{s, h}}(v, \flux; \gamma, \alpha_i) \nonumber\\[5pt]
& := \gamma \, \Big((1 + \alpha_1)\,\| \R_{\rm d} \|_{Q}^2 
      + (1 + \tfrac{1}{\alpha_1})\, \CFriedrichs^2 \, \|\R_{\rm eq}\|^2_{Q} 
			+ \delta_{s, h}\, 
			\big((1 + \alpha_2) \, \| {\dvrg}_x \R_{\rm d}\|^2_{Q} 
			+ (1 + \tfrac{1}{\alpha_2})\,\| \R_{\rm eq} \|^2_{Q}
			\big) \Big),
\label{eq:estimate-deltah-1-non-moving}
\end{alignat}
where $\R_{\rm d}$ and $\R_{\rm eq}$ are defined in \eqref{eq:residuals}, 
$\delta_{s, h}$ is a parameter defined in \eqref{eq:test-function}, 
$\gamma \in \big[\tfrac{1}{2}, +\infty)$, and $\alpha_i >0, i = 1, 2.$
%
A useful particular form of \eqref{eq:estimate-deltah-1-non-moving} arises if we set $\gamma = 1$. Then, 
the estimate has the form
\begin{alignat}{2}
\,\| \nabla_x {e}\|^2_{Q} & + \delta_{s, h} \, \| \partial_t e \|^2_{Q} 
+ \| e \|^2_{\Sigma_T} + \delta_{s, h} \, \| \nabla_x e \|^2_{\Sigma_T} 
\leq \overline{\rm M}^{\rm I}_{\delta_{s, h}}(v, \flux; \alpha_i) \nonumber \\[5pt]
& := (1 + \alpha_1)\,\| \R_{\rm d} \|_{Q}^2 
      + (1 + \tfrac{1}{\alpha_1})\, \CFriedrichs^2 \, \|\R_{\rm eq}\|^2_{Q} 
			+ \delta_{s, h}\, 
			\big( (1 + \alpha_2)\, \| \dvrg_x \R_{\rm d}\|^2_{Q} + 
			        (1 + \tfrac{1}{\alpha_2}) \,\| \R_{\rm eq} \|^2_{Q}
			\big), 
\label{eq:estimate-1-gamma-1}
\end{alignat}
where the best $\alpha_1$ and $\alpha_2$ are defined by relations 
$\alpha^*_1 = \tfrac{\CFriedrichs \, \|\R_{\rm eq}\|_{Q}}{\| \R_{\rm d} \|_{Q}}$ and 
$\alpha^*_2 = \tfrac{\|\R_{\rm eq}\|_{Q}}{\| \dvrg_x \R_{\rm d} \|_{Q}}$.

\end{corollary}
\begin{remark}
In general, $\alpha_1$ and $\alpha_2$ can be positive functions of $t$. Then,  $\alpha^*_1$ and 
$\alpha^*_2$ are also functionals of $t$. In this case, the overall value of the majorant is minimal.
\end{remark}

\begin{corollary}
\label{cor:majorant-2}
Let $v \in V^{\nabla_x \partial_t}_{0, \underline{0}}$
and  ${\flux} \in H^{\dvrg_x, 1}(Q)$. Then, 
%
we have the estimate
\begin{alignat}{2}
\!\!\!
(2 & - \tfrac{1}{\zeta}) (\,\| \nabla_x {e}\|^2_{Q} + \delta_{s, h} \, \| \partial_t e \|^2_{Q}) 
+ \delta_{s, h} \, (1 - \tfrac{1}{\epsilon}) \| \nabla_x e \|^2_{\Sigma_T} + \| e \|^2_{\Sigma_T} 
\leq \overline{\rm M}^{\rm II}_{\delta_{s, h}}(v, \flux; \zeta, \beta_i, \epsilon) \nonumber\\[5pt]
& := \epsilon \, \delta_{s, h} \, \| \R_{\rm d} \|^2_{\Sigma_T} 
+ \zeta \, \Big((1 + \beta_1)\big(
        (1 + \beta_2) \, \| \R_{\rm d} \|^2_{Q} 
				+ (1 + \tfrac{1}{\beta_2}) \CFriedrichs^2 \, \|\R_{\rm eq}\|^2_{Q} \big) 
				+ (1 + \tfrac{1}{\beta_1}) \delta_{s, h}^2 \, \| \partial_t \R_{\rm d} \|^2_Q
      + \delta_{s, h} \, \| \R_{\rm eq} \|_{Q}^2 \Big),
\label{eq:estimate-2-non-moving-domain}
\end{alignat}
where $\R_{\rm d}$ and $\R_{\rm eq}$ are defined in \eqref{eq:residuals}, 
$\delta_{s, h}$ is a parameter defined in \eqref{eq:test-function},  
$\zeta \in \big[\tfrac{1}{2}, +\infty)$, $\epsilon \in [1, +\infty)$, 
and $\beta_i >0, i = 1, 2$. In particular, for $\zeta = 1$ and $\epsilon = 2$, we obtain
\begin{alignat}{2}
\,\| \nabla_x {e}\|^2_{Q} & + \delta_{s, h}\, \| \partial_t e \|^2_{Q} +
\| e \|^2_{\Sigma_T} + \tfrac{\delta_{s, h}}{2} \, \| \nabla_x e \|^2_{\Sigma_T} 
\leq \overline{\rm M}^{\rm II}_{\delta_{s, h}}(v, \flux; \beta_i) \nonumber \\[5pt]
& := 2 \, \delta_{s, h} \, \| \R_{\rm d} \|^2_{\Sigma_T} 
        + (1 + \beta_1)\big(
        (1 + \beta_2) \, \| \R_{\rm d} \|^2_{Q} 
				+ (1 + \tfrac{1}{\beta_2}) \CFriedrichs^2 \, \|\R_{\rm eq}\|^2_{Q} \big) 
				+ (1 + \tfrac{1}{\beta_1}) \delta_{s, h}^2 \, \| \partial_t \R_{\rm d} \|^2_Q
      + \delta_{s, h} \, \| \R_{\rm eq} \|_{Q}^2, 
\label{eq:estimate-2}
\end{alignat}
where the optimal parameters are given by relations 
$$
\beta^*_1 = \tfrac{\delta_{s, h} \, \| \partial \R_{\rm d} \|_{Q}}{\sqrt{(1 + \beta^*_2) \, \| \R_{\rm d} \|^2_{Q} 
				+ \big(1 + \tfrac{1}{\beta^*_2}\big) \CFriedrichs^2 \, \|\R_{\rm eq}\|^2_{Q}}}
				\quad {and} \quad 
				\beta^*_2 = \tfrac{\CFriedrichs \| \R_{\rm eq}\|_{Q}}{\| \R_{\rm d} \|_{Q}}.$$

\end{corollary}

\begin{theorem}
Functionals $\overline{\rm M}^{\rm I}$ and  $\overline{\rm M}^{\rm II}$ 
($\overline{\rm M}^{\rm I}_{\delta_{s, h}}$ and $\overline{\rm M}^{\rm II}_{\delta_{s, h}}$) vanish 
if and only if the approximations $v$ and $\flux$ coincide with the exact solution of the problem 
and its exact flux, i.e., $v = u$ and $\flux = \nabla_x u$. 
\end{theorem}
\ProofBegin
The proof is done for $\overline{\rm M}^{\rm I}$ (for functionals $\overline{\rm M}^{\rm II}$, 
$\overline{\rm M}^{\rm I}_{\delta_{s, h}}$, and $\overline{\rm M}^{\rm II}_{\delta_{s, h}}$, it is done 
analogously). In order to prove the existence of
$(v, \flux) \in V_{0, \underline{0}} \times H^{\dvrg_x, 1}(Q)$, minimising 
$\overline{\rm M}^{\rm I}$, it is enough to construct such a pair explicitly, i.e., set $v = u$ and 
$\flux = \nabla_x u$. In this case, we have $\R_{\rm eq}(v, \flux) = 0$ and $\R_{\rm d}(v, \flux) = 0$. 
Since majorant $\overline{\rm M}^{\rm I}$ is nonnegative, this choice of $v$ and $\flux$ corresponds
to the minimiser. 
On the other hand, if $\overline{\rm M}^{\rm I} = 0$, then $\R_{\rm eq} = 0$ and $\R_{\rm d} = 0$, 
from which it follows that $v$ is the solution of the problem 
\eqref{eq:equation}--\eqref{eq:initial-condition}. By using the argument of uniqueness 
of the solution of the I-BVP, we find that $v = u$. Then, from $\R_{\rm d} = 0$, it automatically follows 
that $\flux = \nabla_x u$.
\ProofEnd
\vskip 10pt
%
\begin{remark}
For the case $\mu = 0$, the majorants presented in Theorems 
\ref{th:theorem-majorant-general-1} and \ref{th:theorem-majorant-general-2}
coincide with the estimates derived in \cite{LMR:Repin:2002}. Computational properties of these 
estimates has been studied in \cite{LMR:MatculevichRepin:2014} and \cite{LMR:Matculevich:2015}.
The paper \cite{LMR:MatculevichRepin:2014} includes two benchmark examples, where error majorants
were applied to approximations reconstructed by the space-time method. 
Numerical results, presented in these examples, confirm the efficient performance of the majorant
(ratios between majorant and error are close to $1$). 
\end{remark}

\section{Modification of majorans $\overline{\rm M}^{\rm I}$ and $\overline{\rm M}^{\rm II}$}
\label{sec:equivalence}

In this section, we deduce modified forms of $\overline{\rm M}^{\rm I}$ and $\overline{\rm M}^{\rm II}$,
which, in general, provide sharper bounds of the error. These estimates contain an additional 
`free' function 
$w \in V^{\Delta_x}_{0, \underline{0}}$. First, we rewrite 
\eqref{eq:energy-balance-equation-with-flux} as follows:
\begin{multline}
	\; \;
	\lambda \, \| \nabla_x{e} \|^2_Q 
	+ \mu \, \| \partial_t \, e \|^2_Q
	+ \tfrac12 \, ( \mu \, \| \nabla_x{e} \|^2_{\Sigma_T} + \lambda \| e \|^2_{\Sigma_T}) \\
	\qquad \;
	= \lambda \, \big( (\widetilde{\R}_{\rm eq}(v; \flux, w) ,  e)_{Q} 
	                           + (\widetilde{\R}_{\rm d}(v; \flux, w) , \nabla_x {e})_{Q}\big) 	
	    + \mu \, \big( (\widetilde{\R}_{\rm eq}(v; \flux, w) , \partial_t {e} )_Q 
	                          + (\widetilde{\R}_{\rm d}(v; \flux, w), \nabla_x \partial_t {e})_{Q} 
	                 \big) \\
	    + \lambda \, \big( (\partial_t w,  e)_Q - (\nabla_x w, \nabla_x {e})_{Q} \big) 
	    + \mu \, \big( (\partial_t w, \partial_t {e})_Q - (\nabla_x w, \nabla_x \partial_t {e})_{Q} \big),
	    \qquad
	    \label{eq:energy-balance-equation-with-flux-with-varphi}
\end{multline}
%
where 
\begin{alignat*}{2}
	\widetilde{\R}_{\rm eq}(v; \flux, w) := f + \dvrg_x \flux - \partial_t (v + w) \quad \mbox{and} \quad 
	\widetilde{\R}_{\rm d}(v; \flux, w) := \flux - \nabla_x (v - w)
\end{alignat*}
are modified residuals. All deducted terms are compensated by adding 
$$\mathcal{I} (e; y, w) 
:= (\partial_t w,  \lambda \, e + \mu \, \partial_t {e})_Q 
- (\nabla_x w, \nabla_x (\lambda \, e + \mu \, \partial_t {e}))_Q.$$
It is not difficult to see that 
\begin{alignat}{2}
\mathcal{I} (e; w) 
& := (\partial_t w,  \lambda \, e + \mu \, \partial_t {e})_Q  
- (\nabla_x w, \nabla_x (\lambda \, e + \mu \, \partial_t {e}))_Q \nonumber\\
& = \lambda \, \big((w \, \vectorn_t, e)_{\Sigma_T} - (w, \partial_t e)_Q \big)
      + \mu \, (\partial_t w, \partial_t {e})_Q
      - \lambda \, (\nabla_x w,  \nabla_x e)_Q
      - \mu \, \big( (\nabla_x w, \nabla_x e \, \vectorn_t)_{\Sigma_T} 
      - (\nabla_x \partial_t w, \nabla_x e)_Q\big) \nonumber\\
& = \lambda \, (w, e)_{\Sigma_T} - \mu \, (\nabla_x w, \nabla_x e)_{\Sigma_T}
- \lambda \, \big((w, \partial_t \, e)_Q + (\nabla_x w,  \nabla_x e)_Q \big)
+ \mu \,\big((\partial_t w, \partial_t {e})_Q + (\nabla_x \partial_t w, \nabla_x e)_Q \big) \nonumber\\
& = \lambda \, (w, e)_{\Sigma_T} - \mu \, (\nabla_x w, \nabla_x e)_{\Sigma_T}
+ \lambda \, \big((w, \partial_t v)_Q + (\nabla_x w,  \nabla_x v)_Q \big)
- \mu \,\big ( (\partial_t w, \partial_t {v})_Q + (\nabla_x \partial_t w, \nabla_x v )_Q \big) \nonumber\\
& \qquad \qquad \qquad \qquad \qquad \qquad \qquad 
- \lambda \, \big((w, \partial_t u)_Q + (\nabla_x w,  \nabla_x u)_Q \big)
+ \mu \,\big (\partial_t w, \partial_t {u})_Q + (\nabla_x \partial_t w, \nabla_x u)_Q \big) \nonumber\\
& = \lambda \, (w, e)_{\Sigma_T} - \mu \, (\nabla_x w, \nabla_x e)_{\Sigma_T}
+ (\lambda w - \mu \partial_t w , \partial_t v)_Q 
+ \big (\nabla_x (\lambda w - \mu \partial_t w),  \nabla_x v \big)_Q 
- (\lambda \, w - \mu \, \partial_t w , f)_Q \nonumber\\
& = \lambda \, (w, e)_{\Sigma_T} - \mu \, (\nabla_x w, \nabla_x e)_{\Sigma_T}
+ \mathcal{J}(v, w)
\label{eq:remaining-terms}
\end{alignat}
%
Here, $$\mathcal{J}(v, w) := (\lambda w - \mu \partial_t w , \partial_t v)_Q 
+ \big (\nabla_x (\lambda w - \mu \partial_t w),  \nabla_x v \big)_Q 
- (\lambda \, w - \mu \, \partial_t w , f)_Q$$ 
is a linear functional that mimics the residual of 
\eqref{eq:variational-formulation} with the test function $\lambda \, w - \mu \, \partial_t w$. The first two 
terms in the RHS of \eqref{eq:remaining-terms} can be estimated as follows:
\begin{equation}
\lambda (w, e)_{\Sigma_T} 
\leq \tfrac{\lambda}{2} \, (\rho_1 \, \| w\|^2_{\Sigma_T} + \tfrac{1}{\rho_1} \| e\|^2_{\Sigma_T})
\label{eq:young-1}
\end{equation}
and 
\begin{equation}
\mu \, (\nabla_x w, \nabla_x e)_{\Sigma_T}
\leq \tfrac{\mu}{2} \,(\rho_2 \, \| \nabla_x w \|^2_{\Sigma_T} 
+ \tfrac{1}{\rho_2} \| \nabla_x e \|^2_{\Sigma_T}),
\label{eq:young-2}
\end{equation}
where $\rho_1, \rho_2 >0$. From \eqref{eq:energy-balance-equation-with-flux-with-varphi}, 
\eqref{eq:remaining-terms}, \eqref{eq:young-1}, and \eqref{eq:young-2}, we obtain
\begin{alignat}{2}
 	\lambda \, \| \nabla_x{e} \|^2_Q 
	& + \mu \, \| \partial_t \, e \|^2_Q
	+ \tfrac12 \, ( \mu \, (1 - \tfrac{1}{\rho_2}) \, \| \nabla_x{e} \|^2_{\Sigma_T} 
	                    + \lambda\, (1 - \tfrac{1}{\rho_1}) \| e \|^2_{\Sigma_T}) 
        \leq \tfrac{\rho_1}{2} \, \| w\|^2_{\Sigma_T}  
	+ \tfrac{\rho_2}{2} \, \| \nabla_x w \|^2_{\Sigma_T} \nonumber\\
	& + \lambda \, \big(
                  ( \widetilde{\R}_{\rm eq}(y, w) ,  e)_Q 
	      + (\widetilde{\R}_{\rm d}(y, w) , \nabla_x {e})_{Q} 
	      \big) 	
	      + \mu \, \big( 
	      (\widetilde{\R}_{\rm eq}(y, w), \partial_t {e} )_Q 
	      + \widetilde{\R}_{\rm d}(y, w), \nabla_x \partial_t {e})_{Q} 
	      \big) + \mathcal{J}(v, w).
\label{eq:energy-balance-inequality-with-flux-with-varphi}
\end{alignat}
By means of this inequality and arguments analogous those used in the proofs of Theorems  
\ref{th:theorem-majorant-general-1} and \ref{th:theorem-majorant-general-2}, 
we deduce advanced forms of majorants.

\begin{theorem}
\label{th:theorem-majorants-second-form} 

(i) For any 
$v, w \in V^{\Delta_x}_{0, \underline{0}}
$ 
and any ${\flux} \in H^{\dvrg_x, 0}(Q)$, 
the alternative estimate holds:
\begin{alignat}{2}
\!\!\!
(2- \tfrac{1}{\gamma}) & 
  (\lambda\,\| \nabla_x e \|^2_{Q} + \mu  \, \| \partial_t e \|^2_{Q}) 
+ \mu \, (1 - \tfrac{1}{\rho_2}) \, \| \nabla_x{e} \|^2_{\Sigma_T} 
      + \lambda\, (1 - \tfrac{1}{\rho_1}) \| e \|^2_{\Sigma_T} \nonumber \\[5pt]
& =: |\!|\!| e  |\!|\!|^2_{\rm I} \leq \overline{\rm M}^{\rm I}_w (v, \flux, w; \gamma, \alpha_i, \rho_i)  := 
         \rho_1 \, \| w\|^2_{\Sigma_T}  + \rho_2 \, \| \nabla_x w \|^2_{\Sigma_T} 
	+ 2\, \mathcal{J}(v, w) \nonumber \\[5pt]
& \qquad + \gamma\Big\{ \lambda \Big((1 + \alpha_1)\,\| \widetilde{\R}_{\rm d} \|_{Q}^2 
      + (1 + \tfrac{1}{\alpha_1})\, \CFriedrichs^2 \, \|\widetilde{\R}_{\rm eq}\|^2_{Q} \Big)
			+ \mu \Big( 
			(1 + \alpha_2) \, \| \dvrg_x \widetilde{\R}_{\rm d}\|^2_{Q} 
			+ (1 + \tfrac{1}{\alpha_2})\,\| \widetilde{\R}_{\rm eq} \|^2_{Q}
			\Big) \Big\}, 
\label{eq:estimate-1-non-moving-h}
\end{alignat}
where $\rho_1, \rho_2, \epsilon \in [1, +\infty)$ and $\alpha_i > 0, i = 1, 2$, are auxiliary parameters.
\vskip 10pt

\noindent
(ii) For any $v, w \in V^{\nabla_x \partial_t}_{0, \underline{0}}$ and any ${\flux} \in H^{\dvrg_x, 1}(Q)$, 
the following inequality holds:
\begin{alignat}{2}
\!\!\!\!
(2 - \tfrac{1}{\zeta}&) (\lambda \,\| \nabla_x {e}\|^2_{Q} + \mu \, \| \partial_t e \|^2_{Q}) 
+ \mu \, (1 - \tfrac{1}{\epsilon} - \tfrac{1}{\rho_2}) \| \nabla_x e \|^2_{\Sigma_T} 
+ \lambda\, (1 - \tfrac{1}{\rho_1}) \,  \| e \|^2_{\Sigma_T} \nonumber\\[5pt]
&  =: |\!|\!| e  |\!|\!|^2_{\rm II} \leq \overline{\rm M}^{\rm II}_w(v, \flux, w; \zeta, \beta_i, \epsilon, \rho_i) 
:= \rho_1 \, \| w\|^2_{\Sigma_T}  
    + \rho_2 \, \| \nabla_x w \|^2_{\Sigma_T} 
    + \epsilon \, \mu \| \widetilde{\R}_{\rm d} \|^2_{\Sigma_T} 
    + 2 \, \mathcal{J}(v, w)
\nonumber\\[5pt]
& \qquad \quad
+ \zeta
  \Big( \lambda \big((1 + \beta_1)\big(
        (1 + \beta_2) \, \| \widetilde{\R}_{\rm d} \|^2_{Q} 
				+ (1 + \tfrac{1}{\beta_2}) \CFriedrichs^2 \, \|\widetilde{\R}_{\rm eq}\|^2_{Q} \big) 
				+ (1 + \tfrac{1}{\beta_1}) \tfrac{\mu^2}{\lambda^2} \, \| \partial_t \widetilde{\R}_{\rm d} \|^2\big)
      + \mu \, \| \R_{\rm eq} \|_{Q}^2 \Big),
\label{eq:estimate-2-non-moving-domain-h}
\end{alignat}
where $\rho_1, \rho_2, \epsilon \in [1, +\infty)$, such that the combination 
$1 - \tfrac{1}{\epsilon} - \tfrac{1}{\rho_2} \geq 0$, $\gamma \in \big[\tfrac{1}{2}, +\infty)$,  
$\zeta \in \big[\tfrac{1}{2}, +\infty)$, and $\beta_i >0, i = 1, 2$.
\noindent
In both inequalities \eqref{eq:estimate-1-non-moving} and \eqref{eq:estimate-2-non-moving-domain}, 
$\widetilde{\R}_{\rm eq}(v; y, w)$ and $\widetilde{\R}_{\rm d}(v; y, w)$ 
are modified residual functionals that follow from \eqref{eq:equation}, and 
$\lambda$ and $\mu$ are positive weights introduced in 
\eqref{eq:upwind-test}. 
\vskip 10pt

\noindent
(iii)
Majorants $\overline{\rm M}^{\rm I}_w$ and $\overline{\rm M}^{\rm II}_w$ satisfy the same properties
that are valid for majorants $\overline{\rm M}^{\rm I}$ and $\overline{\rm M}^{\rm II}$, i.e., they vanish if and only if $v = u$, $\flux = \nabla_x u$, and $w = 0$. 
Moreover, the relation between both forms of the majorants can be written as follows:
\begin{equation}
\overline{\rm M}^{\rm I} = \inf\limits_{w \in V^{\Delta_x}{0, \underline{0}}} 
\overline{\rm M}^{\rm I}_w \quad \mbox{and}
\quad
\overline{\rm M}^{\rm II} = 
\inf\limits_{w \in V^{\partial_t \nabla_x}{0, \underline{0}}} \overline{\rm M}^{\rm II}_w.
\end{equation}

\end{theorem}
%
\ProofBegin
The detailed proofs can be found in works \cite{LMR:Repin:2002, 
LMR:MatculevichNeitaanmakiRepin:2015, LMR:MatculevichRepinPoincare:2014}.
%
\ProofEnd
\vskip 10pt



\noindent
Let us prove the equivalence of modified 
estimate to the error measured in the energy norm \eqref{eq:error-general} and majorant $\overline{\rm M}^{\rm I}_w$. 
Assume that $\flux = \nabla_x u$ and $w = u - v$, then
\begin{alignat}{2}
\!\!\!
(2 & - \tfrac{1}{\gamma}) 
  (\lambda\,\| \nabla_x e \|^2_{Q} + \mu  \, \| \partial_t e \|^2_{Q}) 
+ \mu \, (1 - \tfrac{1}{\rho_2}) \, \| \nabla_x{e} \|^2_{\Sigma_T} 
      + \lambda\, (1 - \tfrac{1}{\rho_1}) \| e \|^2_{\Sigma_T} \nonumber \\[5pt]
& \leq \overline{\rm M}^{\rm I}_w (v, \nabla_x u, u-v; \gamma, \alpha_i, \rho_i)  := 
         \rho_1 \, \| u - v\|^2_{\Sigma_T}  + \rho_2 \, \| \nabla_x (u - v) \|^2_{\Sigma_T} 
	+\mathcal{J}(v, u - v) \nonumber \\[5pt]
& \qquad + \gamma\Big\{ \lambda \Big((1 + \alpha_1)\,\| \nabla_x u - \nabla_x (v - u - v) \|_{Q}^2 
      + (1 + \tfrac{1}{\alpha_1})\, \CFriedrichs^2 \, 
      \| f + \dvrg_x (\nabla_x u) - \partial_t (v + u - v) \|^2_{Q} \Big) \nonumber\\
& \qquad \qquad 
			+ \mu \Big( 
			(1 + \alpha_2) \, \| \dvrg_x (\nabla_x u - \nabla_x (v - u - v))\|^2_{Q} 
			+ (1 + \tfrac{1}{\alpha_2})\,\| f + \dvrg_x (\nabla_x u) - \partial_t (v + u - v) \|^2_{Q} 
			\Big) \Big\} \nonumber\\[5pt]
& = \rho_1 \, \| u - v \|^2_{\Sigma_T}  + \rho_2 \, \| \nabla_x (u - v) \|^2_{\Sigma_T}
   + 2\, \mathcal{J}(v, u - v).
\label{eq:equivalence-2-non-moving}
\end{alignat}
Consider separately the linear functional in the RHS of \eqref{eq:equivalence-2-non-moving}:
\begin{alignat}{2}
\mathcal{J}(v, u - v) & 
=
(\lambda (u - v) - \mu \, \partial_t (u - v) , \partial_t v)_Q 
+ \big (\nabla_x (\lambda (u - v) - \mu \, \partial_t (u - v)),  \nabla_x v \big)_Q 
- (\lambda \, (u - v) - \mu \, \partial_t (u - v) , f)_Q \nonumber\\
& =
(\lambda (u - v) - \mu \, \partial_t (u - v) , \partial_t v)_Q 
+ \big (\nabla_x (\lambda (u - v) - \mu \partial_t (u - v)),  \nabla_x v \big)_Q \nonumber\\
& \qquad \qquad \qquad  \qquad \qquad \qquad  \qquad \quad
- (\lambda \, (u - v) - \mu \, \partial_t (u - v) , \partial_t u)_Q 
- (\lambda \, (u - v) - \mu \, \partial_t (u - v) , \nabla_x u)_Q \nonumber\\
& =
(\lambda (u - v) - \mu \, \partial_t (u - v) ,  - \partial_t (u - v))_Q 
+ \big (\nabla_x (\lambda (u - v) - \mu \partial_t (u - v)),  \nabla_x (u - v) \big)_Q \nonumber\\
& = \lambda \| u - v \|^2_{\Sigma_T}
+  \mu \, \| \partial_t (u - v) \|^2_Q 
+ \lambda \, \| \nabla_x (u - v)\|^2_Q
+ \mu \| \partial_t (u - v)\|^2_{\Sigma_T}.
\label{eq:transformation-of-J-uv}
\end{alignat}
%
In view of \eqref{eq:transformation-of-J-uv} and 
identity \eqref{eq:equivalence-2-non-moving}, we obtain 
\begin{alignat}{2}
\!\!\!
 |\!|\!| e  & |\!|\!|^2_{\rm I} := (2- \tfrac{1}{\gamma}) 
  (\lambda\,\| \nabla_x e \|^2_{Q} + \mu  \, \| \partial_t e \|^2_{Q}) 
+ \mu \, (1 - \tfrac{1}{\rho_2}) \, \| \nabla_x{e} \|^2_{\Sigma_T} 
      + \lambda\, (1 - \tfrac{1}{\rho_1}) \| e \|^2_{\Sigma_T}\nonumber \\[5pt]
& \leq \overline{\rm M}^{\rm I}_w (v, \nabla_x u, e)  := (2\, \lambda + \rho_1)\, \| e \|^2_{\Sigma_T}  
+(2\, \mu + \rho_2)\, \| \nabla_x e \|^2_{\Sigma_T}
+  2\, (\mu \, \| \partial_t e \|^2_Q + \lambda \, \| \nabla_x e\|^2_Q) \nonumber\\[5pt]
& \leq \max \Big\{ \tfrac{2\, \gamma}{2\, \gamma - 1}, 
                  \tfrac{\rho_1 \, (2 \, \lambda + {\rho_1})}{\lambda \, (\rho_1 - 1)}, 
                  \tfrac{\rho_2 \, (2 \, {\mu} + {\rho_2})}{{\mu} \, (\rho_2 - 1)}
                \Big\} 
\Big( (2 - \tfrac{1}{\gamma}) 
  (\lambda\,\| \nabla_x e \|^2_{Q} + \mu  \, \| \partial_t e \|^2_{Q}) 
+ \mu \, (1 - \tfrac{1}{\rho_2}) \, \| \nabla_x{e} \|^2_{\Sigma_T} 
      + \lambda\, (1 - \tfrac{1}{\rho_1}) \| e \|^2_{\Sigma_T} \Big) \nonumber\\[5pt]
& = \max \Big\{ \tfrac{2\, \gamma}{2\, \gamma - 1}, 
                  \tfrac{\rho_1 \, (2 \, \lambda + {\rho_1})}{\lambda \, (\rho_1 - 1)}, 
                  \tfrac{\rho_2 \, (2 \, {\mu} + {\rho_2})}{{\mu} \, (\rho_2 - 1)}
                \Big\} \, |\!|\!| e  |\!|\!|^2_{\rm I}.
\label{eq:equivalence-non-moving}
\end{alignat}
%
The double inequality \eqref{eq:equivalence-non-moving} states the equivalence of 
$|\!|\!| e  |\!|\!|^2_{\rm I}$ and 
$\inf\limits_{\substack{\flux \in H^{\dvrg_x, 0}(Q), \\ w \in V^{\Delta_x}_{0, \underline{0}}
}}
\overline{\rm M}^{\rm I}_w (v, \flux, w; \gamma, \alpha_i, \rho_i)$, i.e., 
\begin{equation}
1 
\leq 
{
\inf\limits_{\substack{\flux \in H^{\dvrg_x, 0}(Q), \\ w \in 
V^{\Delta_x}_{0, \underline{0}}}}
\overline{\rm M}^{\rm I}_w (v, \flux, w; \gamma, \alpha_i, \rho_i)} \Big/ {|\!|\!| e  |\!|\!|^2_{\rm I}}
\leq 
C^{\rm I}_{\rm eq},
\end{equation}
where the constant $C^{\rm I}_{\rm eq} := \max \Big\{ \tfrac{2\, \gamma}{2\, \gamma - 1}, 
                  \tfrac{\rho_1 \, (2 \, \lambda + {\rho_1})}{\lambda \, (\rho_1 - 1)}, 
                  \tfrac{\rho_2 \, (2 \, {\mu} + {\rho_2})}{{\mu} \, (\rho_2 - 1)}
                \Big\}$ is explicitly computable. 

Analogous equivalence of $|\!|\!| e  |\!|\!|^2_{\rm II}$ and majorant
$\inf\limits_{\substack{\flux \in H^{\dvrg_x, 1}(Q), \\ w \in V^{\nabla_x \partial_t}_{0, \underline{0}}
}}
\overline{\rm M}^{\rm II}_w (v, \flux, w; \gamma, \beta_i, \epsilon, \rho_i)$ can 
be formulated as follows:
\begin{equation}
1 
\leq 
{
\inf\limits_{\substack{\flux \in H^{\dvrg_x, 1}(Q), \\ w \in V^{\nabla_x \partial_t}_{0, \underline{0}}
}}
\overline{\rm M}^{\rm II}_w (v, \flux, w; \gamma, \beta_i, \epsilon, \rho_i)} \Big/ {|\!|\!| e  |\!|\!|^2_{\rm II}}
\leq 
C^{\rm II}_{\rm eq},
\end{equation}
where $C^{\rm II}_{\rm eq} := \max \Big\{ \tfrac{2\, \gamma}{2\, \gamma - 1}, 
                  \tfrac{\rho_1 \, (2\, \lambda + {\rho_1})}{{\lambda} \, (\rho_1 - 1)}, 
                  \tfrac{(2 \, {\mu} + {\rho_2})}{{\mu} \, \big(1 - \tfrac{1}{\epsilon} - \tfrac{1}{\rho_2}\big)}
                \Big\}$ is explicitly defined.

%

\section{Acknowledgements}

The research is supported by the Austrian Science Fund (FWF) through the NFN S117-03 project. 
The third author of the paper, Prof. Sergey Repin, expresses deep gratitude to RICAM Institute for 
the kind hospitality during the Special Semester. 

\bibliographystyle{plain}
\bibliography{bib/lib}

\begin{thebibliography}{10}

\bibitem{AinsworthOden2000}
M.~Ainsworth and J.~T. Oden.
\newblock {\em A posteriori error estimation in finite element analysis}.
\newblock Wiley and Sons, New York, 2000.

\bibitem{LMR:BabuskaJanik:1989}
I.~Babu{\v{s}}ka and T.~Janik.
\newblock The {$h$}-{$p$} version of the finite element method for parabolic
  equations. {I}. {T}he {$p$}-version in time.
\newblock {\em Numer. Methods Partial Differential Equations}, 5(4):363--399,
  1989.

\bibitem{LMR:BabuskaJanik:1990}
I.~Babu{\v{s}}ka and T.~Janik.
\newblock The {$h$}-{$p$} version of the finite element method for parabolic
  equations. {II}. {T}he {$h$}-{$p$} version in time.
\newblock {\em Numer. Methods Partial Differential Equations}, 6(4):343--369,
  1990.

\bibitem{LMR:BankSmith:1993}
R.E. Bank and R.K. Smith.
\newblock A posteriori error estimates based on hierarchical bases.
\newblock {\em SIAM J. Numer. Anal.}, 30(4):921--935, 1993.

\bibitem{LMR:Bazilevsetal2006}
Y.~Bazilevs, L.~Beir{\~a}o~da Veiga, J.~A. Cottrell, T.~J.~R. Hughes, and
  G.~Sangalli.
\newblock Isogeometric analysis: approximation, stability and error estimates
  for {$h$}-refined meshes.
\newblock {\em Math. Models Methods Appl. Sci.}, 16(7):1031--1090, 2006.

\bibitem{LMR:Bazilevsetall2010}
Y.~Bazilevs, V.~M. Calo, J.~A. Cottrell, J.~A. Evans, T.~J.~R. Hughes,
  S.~Lipton, M.~A. Scott, and T.~W. Sederberg.
\newblock Isogeometric analysis using {T}-splines.
\newblock {\em Comput. Methods Appl. Mech. Engrg.}, 199(5-8):229--263, 2010.

\bibitem{LMR:Veigaetall2011}
L.~Beir{\~a}o~da Veiga, A.~Buffa, D.~Cho, and G.~Sangalli.
\newblock Iso{G}eometric analysis using {T}-splines on two-patch geometries.
\newblock {\em Comput. Methods Appl. Mech. Engrg.}, 200(21-22):1787--1803,
  2011.

\bibitem{LMR:BeiraodaVeigaBuffaRivasSangalli:2011}
L.~Beir{\~a}o~da Veiga, A.~Buffa, J.~Rivas, and G.~Sangalli.
\newblock Some estimates for {$h$}-{$p$}-{$k$}-refinement in isogeometric
  analysis.
\newblock {\em Numer. Math.}, 118(2):271--305, 2011.

\bibitem{LMR:Bressan2013}
A.~Bressan.
\newblock Some properties of {LR}-splines.
\newblock {\em Comput. Aided Geom. Design}, 30(8):778--794, 2013.

\bibitem{LMR:BuffaGiannelli2015}
A.~Buffa and C.~Giannelli.
\newblock Adaptive isogeometric methods with hierarchical splines: error
  estimator and convergence.
\newblock Technical Report arxiv: 1502.00565, arxiv:math.NA, 2015.

\bibitem{LMR:DedeSantos2012}
L.~Ded{\`e} and H.~A. F.~A. Santos.
\newblock B-spline goal-oriented error estimators for geometrically nonlinear
  rods.
\newblock {\em Comput. Mech.}, 49(1):35--52, 2012.

\bibitem{LMR:Deshpandeetall:1995}
A.~Deshpande, S.~Malhotra, M.H. Schultz, and C.C. Douglas.
\newblock A rigorous analysis of time domain parallelism.
\newblock {\em Parallel Algorithms and Applications}, 6(1):53--62, 1995.

\bibitem{LMR:DokkenLychePettersen2013}
T.~Dokken, T.~Lyche, and K.~F. Pettersen.
\newblock Polynomial splines over locally refined box-partitions.
\newblock {\em Comput. Aided Geom. Design}, 30(3):331--356, 2013.

\bibitem{LMR:DorfelJuttlerSimeon2010}
M.~R. D{\"o}rfel, B.~J{\"u}ttler, and B.~Simeon.
\newblock Adaptive isogeometric analysis by local {$h$}-refinement with
  {T}-splines.
\newblock {\em Comput. Methods Appl. Mech. Engrg.}, 199(5-8):264--275, 2010.

\bibitem{LMR:EvansHughes2013}
J.~A. Evans and T.~J.~R. Hughes.
\newblock Explicit trace inequalities for isogeometric analysis and parametric
  hexahedral finite elements.
\newblock {\em Numer. Math.}, 123(2):259--290, 2013.

\bibitem{LMR:GaevskayaRepin:2005}
A.~V. Gaevskaya and S.~I. Repin.
\newblock A posteriori error estimates for approximate solutions of linear
  parabolic problems.
\newblock {\em Springer, Differential Equations}, 41(7):970--983, 2005.

\bibitem{LMR:Gander:2015}
M.~Gander.
\newblock 50 years of time parallel time integration.
\newblock In {\em Multiple Shooting and Time Domain Decomposition}, volume~16,
  pages 69--114. Springer-Verlag, Berlin, 2015.
\newblock Theory, algorithm, and applications.

\bibitem{LMR:GanderNeumuller:2014}
M.~J. Gander and M.~Neum{\"u}ller.
\newblock Analysis of a time multigrid algorithm for dg-discretizations in
  time.
\newblock Technical Report NuMa-Report 2014-07, Johannes Kepler University
  Linz, Institute for Computational Mathematics, Linz, Linz, 2014.

\bibitem{LMR:GanderNeumuller:2016}
M.~J. Gander and M.~Neum{\"u}ller.
\newblock Analysis of a new space-time parallel multigrid algorithm for
  parabolic problems.
\newblock {\em SIAM J. Sci. Comput.}, 38(4):A2173--A2208, 2016.

\bibitem{LMR:GiannelliJuttlerSpeleers2012}
C.~Giannelli, B.~J{\"u}ttler, and H.~Speleers.
\newblock T{HB}-splines: the truncated basis for hierarchical splines.
\newblock {\em Comput. Aided Geom. Design}, 29(7):485--498, 2012.

\bibitem{LMR:Hackbusch:1984}
W.~Hackbusch.
\newblock Parabolic multigrid methods.
\newblock In {\em Computing methods in applied sciences and engineering, {VI}
  ({V}ersailles, 1983)}, pages 189--197. North-Holland, Amsterdam, 1984.

\bibitem{LMR:Hansbo:1994}
P.~Hansbo.
\newblock Space-time oriented streamline diffusion methods for nonlinear
  conservation laws in one dimension.
\newblock {\em Comm. Numer. Meth. Eng.}, 10(3):203--215, 1994.

\bibitem{LMR:Johannessen:2009}
K.~A. Johannessen.
\newblock An adaptive isogeometric finite element analysis.
\newblock Technical report, Master Thesis, Norwegian University of Science and
  Technology, 2009.

\bibitem{LMR:Johnson:1987}
C.~Johnson.
\newblock {\em Numerical solution of partial differential equations by the
  finite element method}.
\newblock Dover Publications, Inc., Mineola, NY, 1987.

\bibitem{LMR:JohnsonSaranen:1986}
C.~Johnson and J.~Saranen.
\newblock Streamline diffusion methods for the incompressible euler and
  navier-stokes equations.
\newblock {\em Math. Comp.}, 47(175):1--18, 1986.

\bibitem{LMR:Karabelas:2015}
E.~Karabelas.
\newblock {\em Space-time discontinuous Galerkin methods for cardic
  electro-mechanics}.
\newblock PhD thesis, Technische Universitat Graz, 2015.

\bibitem{LMR:KarabelasNeumuller:2015}
E.~Karabelas and M.~Neum{\"u}ller.
\newblock Generating admissible space-time meshes for moving domains in
  $(d+1)$-dimensions.
\newblock Technical Report NuMa-Report 2015-07, Johannes Kepler University
  Linz, Institute for Computational Mathematics, Linz, Linz, 2015.

\bibitem{LMR:KleissJuttlerZulehner:2012}
S.~K. Kleiss, B.~Jüttler, and W.~Zulehner.
\newblock Enhancing isogeometric analysis by a finite element?based local
  refinement strategy.
\newblock {\em Comput. Methods Appl. Mech. Engrg.}, 213--216:168--182, 2012.

\bibitem{LMR:KleissTomar2015}
S.~K. Kleiss and S.~K. Tomar.
\newblock Guaranteed and sharp a posteriori error estimates in isogeometric
  analysis.
\newblock {\em Comput. Math. Appl.}, 70(3):167--190, 2015.

\bibitem{LMR:Kraft1997}
R.~Kraft.
\newblock Adaptive and linearly independent multilevel {$B$}-splines.
\newblock In {\em Surface fitting and multiresolution methods
  ({C}hamonix--{M}ont-{B}lanc, 1996)}, pages 209--218. Vanderbilt Univ. Press,
  Nashville, TN, 1997.

\bibitem{LMR:KumarKvamsdalJohannessen2015}
M.~Kumar, T.~Kvamsdal, and K.~A. Johannessen.
\newblock Simple a posteriori error estimators in adaptive isogeometric
  analysis.
\newblock {\em Comput. Math. Appl.}, 70(7):1555--1582, 2015.

\bibitem{LMR:Kuru:2013}
G.~Kuru.
\newblock Goal-adaptive isogeometric analysis with hierarchical splines.
\newblock Technical report, Master's thesis, Mechanical Engineering, Eindhoven
  University of Technology, 2013.

\bibitem{LMR:Kuruetall2014}
G.~Kuru, C.~V. Verhoosel, K.~G. van~der Zee, and E.~H. van Brummelen.
\newblock Goal-adaptive isogeometric analysis with hierarchical splines.
\newblock {\em Comput. Methods Appl. Mech. Engrg.}, 270:270--292, 2014.

\bibitem{LMR:Ladyzhenskaya:1954}
O.~A. Ladyzhenskaya.
\newblock On solvability of classical boundary value problems for equations of
  parabolic and hyperbolic types.
\newblock {\em Dokl. Akad. Nauk SSSR}, 97(3):395--398, 1954.

\bibitem{LMR:Ladyzhenskaya:1985}
O.~A. Ladyzhenskaya.
\newblock {\em The boundary value problems of mathematical physics}.
\newblock Springer, New York, 1985.

\bibitem{LMR:Ladyzhenskayaetal:1967}
O.~A. Ladyzhenskaya, V.~A. Solonnikov, and N.N. Uraltseva.
\newblock {\em Linear and quasilinear equations of parabolic type}.
\newblock Nauka, Moscow, 1967.

\bibitem{LMR:Lang:2001}
J.~Lang.
\newblock {\em Adaptive multilevel solution of nonlinear parabolic {PDE}
  systems}, volume~16 of {\em Lecture Notes in Computational Science and
  Engineering}.
\newblock Springer-Verlag, Berlin, 2001.
\newblock Theory, algorithm, and applications.

\bibitem{LMR:LangerMooreNeumueller:2016a}
U.~Langer, S.~Moore, and M.~Neum\"uller.
\newblock Space-time isogeometric analysis of parabolic evolution equations.
\newblock {\em Comput. Methods Appl. Mech. Engrg.}, 306:342--363, 2016.

\bibitem{LMR:LangerRepinWolfmayr:2015}
U.~Langer, S.~Repin, and M.~Wolfmayr.
\newblock Functional a posteriori error estimates for parabolic time-periodic
  boundary value problems.
\newblock {\em CMAM}, 15(3):353--372, 2015.

\bibitem{LMR:LangerRepinWolfmayr:2016}
U.~Langer, S.~Repin, and M.~Wolfmayr.
\newblock Functional a posteriori error estimates for time-periodic parabolic
  optimal control problems.
\newblock {\em Numer. Func. Anal. Opt.}, 2016.

\bibitem{LMR:LubichOstermann:1987}
Ch. Lubich and A.~Ostermann.
\newblock Multigrid dynamic iteration for parabolic equations.
\newblock {\em BIT}, 27(2):216--234, 1987.

\bibitem{Malietall2014}
O.~Mali, P.~Neittaanm{\"a}ki, and S.~Repin.
\newblock {\em Accuracy verification methods}, volume~32 of {\em Computational
  Methods in Applied Sciences}.
\newblock Springer, Dordrecht, 2014.

\bibitem{LMR:Matculevich:2015}
S.~Matculevich.
\newblock {\em Fully reliable a posteriori error control for evolutionary
  problems}.
\newblock PhD thesis, Jyv{\"a}skyl{\"a} Studies in Computing, University of
  Jyv{\"a}skyl{\"a}, 2015.

\bibitem{LMR:MatculevichNeitaanmakiRepin:2015}
S.~Matculevich, P.~Neittaanm{\"a}ki, and S.~Repin.
\newblock A posteriori error estimates for time-dependent reaction-diffusion
  problems based on the {P}ayne--{W}einberger inequality.
\newblock {\em AIMS}, 35(6):2659--2677, 2015.

\bibitem{LMR:MatculevichRepinPoincare:2014}
S.~Matculevich and S.~Repin.
\newblock Computable bounds of the distance to the exact solution of parabolic
  problems based on {P}oincar\'e type inequalities.
\newblock {\em Zap. Nauchn. Sem. S.-Peterburg. Otdel. Mat. Inst. Steklov
  (POMI)}, 425(1):7--34, 2014.

\bibitem{LMR:MatculevichRepin:2014}
S.~Matculevich and S.~Repin.
\newblock Computable estimates of the distance to the exact solution of the
  evolutionary reaction-diffusion equation.
\newblock {\em Appl. Math. and Comput.}, 247:329--347, 2014.

\bibitem{LMR:MatculevichRepin:2015}
S.~Matculevich and S.~Repin.
\newblock Explicit constants in poincar\'{e}-type inequalities for simplicial
  domains.
\newblock {\em Computational Methods in Applied Mathematics - CMAM},
  16(2):277--298, 2016.

\bibitem{LMR:Mollet:2014}
C.~Mollet.
\newblock Stability of {P}etrov-{G}alerkin discretizations: application to the
  space-time weak formulation for parabolic evolution problems.
\newblock {\em Comput. Methods Appl. Math.}, 14(2):231--255, 2014.

\bibitem{LMR:NeumullerSteinbach:2011}
M.~Neum{\"u}ller and O.~Steinbach.
\newblock Refinement of flexible space-time finite element meshes and
  discontinuous {G}alerkin methods.
\newblock {\em Comput. Vis. Sci.}, 14(5):189--205, 2011.

\bibitem{LMR:NeumullerSteinbach:2013}
M.~Neum{\"u}ller and O.~Steinbach.
\newblock A {DG} space-time domain decomposition method.
\newblock In {\em Domain decomposition methods in science and engineering
  {XX}}, volume~91 of {\em Lect. Notes Comput. Sci. Eng.}, pages 623--630.
  Springer, Heidelberg, 2013.

\bibitem{LMR:NguyenThanhetall:2011}
N.~Nguyen-Thanh, H.~Nguyen-Xuan, S.~P.~A. Bordas, and T.~Rabczuk.
\newblock Isogeometric analysis using polynomial splines over hierarchical
  t-meshes for two-dimensional elastic solids.
\newblock {\em Comput. Methods Appl. Mech. Engrg.}, 200:1892--1908, 2011.

\bibitem{LMR:RepinDeGruyterMonograph:2008}
S.~Repin.
\newblock {\em A posteriori estimates for partial differential equations},
  volume~4 of {\em Radon Series on Computational and Applied Mathematics}.
\newblock Walter de Gruyter GmbH \& Co. KG, Berlin, 2008.

\bibitem{LMR:Repin:2002}
S.~I. Repin.
\newblock Estimates of deviations from exact solutions of initial-boundary
  value problem for the heat equation.
\newblock {\em Rend. Mat. Acc. Lincei}, 13(9):121--133, 2002.

\bibitem{LMR:RepinTomar:2010}
S.~I. Repin and S.~K. Tomar.
\newblock A posteriori error estimates for approximations of evolutionary
  convection-diffusion problems.
\newblock {\em J. Math. Sci. (N. Y.)}, 170(4):554--566, 2010.

\bibitem{LMR:Repin:1997}
S.I. Repin.
\newblock A posteriori error estimation for nonlinear variational problems by
  duality theory.
\newblock {\em Zapiski Nauchnych Seminarov POMIs}, 243:201--214,, 1997.

\bibitem{LMR:Repin:1999}
S.I. Repin.
\newblock A posteriori error estimates for approximate solutions to variational
  problems with strongly convex functionals.
\newblock {\em Journal of Mathematical Sciences}, 97:4311--4328, 1999.

\bibitem{LMR:Repin:2000}
S.I. Repin.
\newblock A posteriori error estimation for variational problems with uniformly
  convex functionals.
\newblock {\em Math. Comput.}, 69(230):481--500, 2000.

\bibitem{LMR:SchwabStevenson:2009}
C.~Schwab and R.~Stevenson.
\newblock Space-time adaptive wavelet methods for parabolic evolution problems.
\newblock {\em Math. Comp.}, 78(267):1293--1318, 2009.

\bibitem{LMR:Scottetall2011}
M.~A. Scott, M.~J. Borden, C.~V. Verhoosel, T.~W. Sederberg, and T.~J.~R.
  Hughes.
\newblock Isogeometric finite element data structures based on {B}\'ezier
  extraction of {T}-splines.
\newblock {\em Internat. J. Numer. Methods Engrg.}, 88(2):126--156, 2011.

\bibitem{LMR:Scottetall2012}
M.~A. Scott, X.~Li, T.~W. Sederberg, and T.~J.~R. Hughes.
\newblock Local refinement of analysis-suitable {T}-splines.
\newblock {\em Comput. Methods Appl. Mech. Engrg.}, 213/216:206--222, 2012.

\bibitem{LMR:Sederbergetall2004}
T.~W. Sederberg, D.~C. Cardon, G.~T. Finnigan, N.~N. North, J.~Zheng, and
  T.~Lyche.
\newblock T-splines simplification and local refinement.
\newblock {\em ACM Trans. Graphics}, 23(3):276--283, 2004.

\bibitem{LMR:Sederbergetall2003}
T.~W. Sederberg, J.~Zheng, A.~Bakenov, and A.~Nasri.
\newblock T-splines and t-nurccs.
\newblock {\em ACM Trans. Graphics}, 22(3):477--484, 2003.

\bibitem{LMR:Steinbach:2015}
O.~Steinbach.
\newblock Space-time finite element methods for parabolic problems.
\newblock {\em Computational Methods in Applied Mathematics}, 15(4):551--566,
  2015.

\bibitem{LMR:TagliabueDedeQuarteroni:2014}
A.~Tagliabue, L.~Ded{\`e}, and A.~Quarteroni.
\newblock Isogeometric analysis and error estimates for high order partial
  differential equations in fluid dynamics.
\newblock {\em Comput. \& Fluids}, 102:277--303, 2014.

\bibitem{LMR:Takizawaetall:2012}
K.~Takizawa, K.~Schjodt, A.~Puntel, N.~Kostov, and T.~E. Tezduyar.
\newblock Patient-specific computer modeling of blood flow in cerebral arteries
  with aneurysm and stent.
\newblock {\em Comput. Mech.}, 50(6):675--686, 2012.

\bibitem{LMR:TakizawaTezduyar:2011}
K.~Takizawa and T.~E. Tezduyar.
\newblock Multiscale space-time fluid-structure interaction techniques.
\newblock {\em Comput. Mech.}, 48(3):247--267, 2011.

\bibitem{LMR:TakizawaTezduyar:2014}
K.~Takizawa and T.~E. Tezduyar.
\newblock Space-time computation techniques with continuous representation in
  time ({ST}-{C}).
\newblock {\em Comput. Mech.}, 53(1):91--99, 2014.

\bibitem{LMR:Thomee:2006}
V.~Thom{\'e}e.
\newblock {\em Galerkin finite element methods for parabolic problems},
  volume~25 of {\em Springer Series in Computational Mathematics}.
\newblock Springer-Verlag, Berlin, second edition, 2006.

\bibitem{LMR:UrbanPatera:2014}
K.~Urban and A.~T. Patera.
\newblock An improved error bound for reduced basis approximation of linear
  parabolic problems.
\newblock {\em Math. Comp.}, 83(288):1599--1615, 2014.

\bibitem{LMR:ZeeVerhoosel2011}
K.~G. van~der Zee and C.~V. Verhoosel.
\newblock Isogeometric analysis-based goal-oriented error estimation for
  free-boundary problems.
\newblock {\em Finite Elem. Anal. Des.}, 47(6):600--609, 2011.

\bibitem{LMR:VandewalleHorton:1995}
S.~Vandewalle and G.~Horton.
\newblock Fourier mode analysis of the multigrid waveform relaxation and
  time-parallel multigrid methods.
\newblock {\em Computing}, 54(4):317--330, 1995.

\bibitem{LMR:Vuongetall2011}
A.-V. Vuong, C.~Giannelli, B.~J{\"u}ttler, and B.~Simeon.
\newblock A hierarchical approach to adaptive local refinement in isogeometric
  analysis.
\newblock {\em Comput. Methods Appl. Mech. Engrg.}, 200(49-52):3554--3567,
  2011.

\bibitem{LMR:Wangetall2011}
P.~Wang, J.~Xu, J.~Deng, and F.~Chen.
\newblock Adaptive isogeometric analysis using rational pht-splines.
\newblock {\em Computer-Aided Design}, 43(11):1438--1448, 2011.

\bibitem{LMR:Zeidler:1990a}
E.~Zeidler.
\newblock {\em Nonlinear functional analysis and its applications. {II}/{A}}.
\newblock Springer-Verlag, New York, 1990.

\end{thebibliography}

\end{document}